\documentclass[11pt]{article}
\usepackage{latexsym}
\usepackage{amsfonts,amssymb,amsmath}

\newtheorem{thm}{Theorem}[section]
\newtheorem{cor}[thm]{Corollary}
\newtheorem{lem}[thm]{Lemma}
\newtheorem{pro}[thm]{Proposition}
\newtheorem{defn}[thm]{Definition}

\title{Infinitely generated semigroups and polynomial complexity} 

\author{ J.C.\ Birget  } 

\date{\footnotesize{\empty}} 
\begin{document}
\maketitle

\begin{abstract}
This paper continues the functional approach to the {\sf P}-versus-{\sf NP}
problem, begun in \cite{s1f}. 
Here we focus on the monoid ${\cal RM}_2^{\sf P}$ of right-ideal morphisms 
of the free monoid, that have polynomial input balance and polynomial 
time-complexity. We construct a machine model for the functions in 
${\cal RM}_2^{\sf P}$, and evaluation functions. We prove that 
${\cal RM}_2^{\sf P}$ is not finitely generated, and use this to show 
separation results for time-complexity.
\end{abstract}

%%%%%%%%%%%%%%%%%%%%%
\section{ Introduction}

In \cite{s1f} we defined the monoids of partial functions {\sf fP} and 
${\cal RM}_2^{\sf P}$. 
The question whether {\sf P} $=$ {\sf NP} is equivalent to the question 
whether these monoids are regular.
The monoid {\sf fP} consists of all partial functions $A^* \to A^*$ that
are computable by deterministic Turing machines in polynomial time, and
that have polynomial I/O-balance. The submonoid ${\cal RM}_2^{\sf P}$ 
consists of the elements of {\sf fP} that are right-ideal morphisms of $A^*$.
One-way functions (according to worst-case time-complexity) are exactly the 
non-regular elements of {\sf fP}. 
It is known that one-way functions (according to worst-case time-complexity) 
exist iff {\sf P} $\neq$ {\sf NP}. 
Also, $f \in {\cal RM}_2^{\sf P}$ is regular in {\sf fP} iff $f$ is regular 
in ${\cal RM}_2^{\sf P}$. 
Hence, {\sf P} $=$ {\sf NP} iff {\sf fP} is regular, iff 
${\cal RM}_2^{\sf P}$ is regular.
We refer to \cite{DuKo,Papadim} for background on {\sf P} and {\sf NP}.

The original motivation for studying ${\cal RM}_2^{\sf P}$ in addition to 
{\sf fP} was that ${\cal RM}_2^{\sf P}$ is reminiscent of the 
Thompson-Higman groups \cite{McKTh,Th,Hig74,CFP,BiThomps,BiCoNP} and the 
Thompson-Higman monoids \cite{JCBmonThH}.
It also quickly turned out that ${\cal RM}_2^{\sf P}$, while having the 
same connection to {\sf P}-vs.-{\sf NP} as {\sf fP}, has different 
properties than {\sf fP} (e.g., regarding the Green relations, and 
actions on $\{0,1\}^{\omega}$; see  \cite{s1f, equiv}). 
It is hard to know whether this approach will contribute to a solution of 
the {\sf P}-vs.-{\sf NP} problem, but the monoids {\sf fP} and 
${\cal RM}_2^{\sf P}$ are interesting by themselves.

\smallskip

Above and in the rest of the paper we use the following notation and 
terminology.
We have an alphabet $A$, which will be $\{0,1\}$ unless the contrary is
explicitly stated, and $A^*$ denotes the set of all strings over $A$, 
including the empty string $\varepsilon$. For $x \in A^*$, $|x|$ denotes the 
length of the string $x$. For a partial function $f: A^* \to A^*$, the 
domain is ${\sf Dom}(f) = \{ x \in A^* : f(x)$ is defined\}, and the image 
is ${\sf Im}(f) = f(A^*) = f({\sf Dom}(f))$. When we say ``function'', we 
mean partial function (except when we explicitly say ``total function''). 
Similarly, for a deterministic input-output Turing machine with input-output
alphabet $A$, the domain of the machine is the set of input words for which 
the machine produces an output; and the set of output words is the image of
the machine.

A function $f: A^* \to A^*$ is called {\em polynomially balanced} iff
there exists polynomials $p, q$ such that for all $x \in {\sf Dom}(f)$:
$|f(x)| \leq p(|x|)$ and $|x| \leq q(|f(x)|)$.
The polynomial $q$ is called an {\em input balance} function for $f$.

As we said already, {\sf fP} is the set of partial functions 
$f: A^* \to A^*$ that are polynomially balanced, and such that
$x \in {\sf Dom}(f) \longmapsto f(x)$ is computable by a deterministic 
polynomial-time Turing machine.
Hence, ${\sf Dom}(f)$ is in {\sf P} when $f \in {\sf fP}$, and it is not
hard to show that ${\sf Im}(f)$ is in {\sf NP}.
Clearly, {\sf fP} is a monoid under function composition.

A function $f: A^* \to A^*$ is said to be {\em one-way (with respect to
worst-case complexity)} iff $f \in {\sf fP}$, but there exists {\em no} 
deterministic polynomial-time algorithm which, on every input 
$y \in {\sf Im}(f)$, outputs some $x \in A^*$ such that $f(x) = y$.  By 
``one-way'' we will always mean one-way with respect to worst-case 
complexity; hence, these functions are not ``cryptographic one-way 
functions'' (in the sense of, e.g., \cite{DiffHell, LevinTale, Goldreich}).
However, they are important for the {\sf P}-vs.-{\sf NP} problem because of
the following folklore fact (see e.g., \cite{HemaOgi} p.\ 33): 
{\em One-way functions exist iff ${\sf P} \neq {\sf NP}$.} 

As is easy to prove (see the Introduction of \cite{s1f}), $f \in {\sf fP}$ 
is not one-way iff $f$ is {\em regular} in {\sf fP}. 
By definition, an element $f$ in a monoid $M$ is regular iff there exists 
$f' \in M$ such that $f f' f = f$; in this case, $f'$ is called an 
{\em inverse} of $f$.\footnote{The terminology varies, depending on the 
field. In semigroup theory $f'$ such that $f f' f = f$ is called a
semi-inverse or a pseudo-inverse of $f$, in numerical mathematics $f'$
is called a generalized inverse, in ring theory and in category theory it's 
called a weak inverse. In semigroup theory the term ``inverse'' of $f$
is only applied to $f'$ if $f' f f' = f'$ holds in addition to $f f' f = f$. 
It is easy to see that if $f f' f = f$ then $f' f f'$ ($= v$) satisfies
$f v f = f$ and $v f v = v$.
}
A monoid $M$ is called regular iff all the elements 
of $M$ are regular.  In summary we have: 
{\em The monoid {\sf fP} is regular iff ${\sf P} = {\sf NP}$.}

\smallskip

Let us look in more detail at the monoid ${\cal RM}_2^{\sf P}$.
A {\em right ideal} of $A^*$ is a subset $R \subseteq A^*$ such that 
$R \, A^* = R$ (i.e., $R$ is closed under right-concatenation by any string).
For two strings $v, w \in A^*$, we say that $v$ is a {\em prefix} of $w$
iff $(\exists x \in A^*)[ \, vx = w \, ]$. 
A {\em prefix code} in $A^*$ is a set $P \subset A^*$ such that no word in 
$P$ is a prefix of another word in $P$.  
For any right ideal $R$ there exists a unique prefix code $P_{_R}$ such that
$R = P_{_{\!\! R}} \, A^*$; we say that $P_{_{\!\! R}}$ generates $R$ as a 
right ideal. For details, see e.g.\ \cite{BiThomps,JCBmonThH}; a good 
reference on prefix codes, and variable-length codes in general is 
\cite{ThCodes}.

A {\em right-ideal morphism} is a partial function $h: A^* \to A^*$ such
that for all $x \in {\sf Dom}(h)$ and all $w \in A^*$: \ $h(xw) = h(x) \, w$.
In that case, ${\sf Dom}(h)$ and ${\sf Im}(h)$ are right ideals. 
For a right-ideal morphism $h$, let ${\sf domC}(h)$ (called 
the {\em domain code}) be the prefix code that generates ${\sf Dom}(h)$ as a  
right ideal.  Similarly, let ${\sf imC}(h)$, called the {\em image code}, be 
the prefix code that generates ${\sf Im}(h)$.  So a right-ideal morphism $h$ 
is determined by $h|_{{\sf domC}(h)}$ (the restriction of $h$ to its domain 
code). In general, ${\sf imC}(h) \subseteq h({\sf domC}(h))$, and it can 
happen that ${\sf imC}(h) \neq h({\sf domC}(h))$.
We define 

\medskip

\hspace{.8in} ${\cal RM}_2^{\sf P} \ = \ \{ f \in {\sf fP} : f$ is a 
                right-ideal morphism of $A^*\}$.

\medskip

\noindent By Prop.\ 2.6 in \cite{s1f}, $f \in {\cal RM}_2^{\sf P}$ is 
regular in ${\cal RM}_2^{\sf P}$ iff $f$ is regular in {\sf fP}. 
{\em The monoid ${\cal RM}_2^{\sf P}$ is regular iff ${\sf P} = {\sf NP}$.}

We saw (Cor.\ 2.9 in \cite{s1f}) that {\sf fP} and ${\cal RM}_2^{\sf P}$ are
not isomorphic, that the group of units of ${\cal RM}_2^{\sf P}$ is 
trivial (Prop.\ 2.12 in \cite{s1f}), and that ${\cal RM}_2^{\sf P}$ has only
one non-0 ${\cal J}$-class (Prop.\ 2.7 in \cite{s1f}). 
In \cite{equiv} we will see that ${\cal RM}_2^{\sf P}$ has interesting 
actions on $\{0,1\}^{\omega}$, and has interesting homomorphic images (some 
of which are regular monoids, and some of which are regular iff {\sf P} $=$
{\sf NP}). Overall, ${\cal RM}_2^{\sf P}$ seems to have ``more structure''
than {\sf fP}.

It is proved in \cite{s1f} (Section 3) that {\sf fP} is isomorphic to a 
submonoid of ${\cal RM}_2^{\sf P}$. To prove this, we use an encoding of 
the three-letter alphabet $\{0,1,\#\}$ into words over the two-letter 
alphabet $\{0,1\}$; this encoding will also be used here. First, we encode 
the alphabet $\{0, 1, \#\}$ by $ \, {\sf code}(0) = 00$,
 \ ${\sf code}(1) = 01$, \ ${\sf code}(\#) = 11$.
A word $x_1 \ldots x_n \in \{0, 1, \#\}^*$ is encoded to 
 \ ${\sf code}(x_1) \ \ldots \ {\sf code}(x_n)$.  
For a fixed $k>0$, a $k$-tuple of words $(u_1, \ldots, u_{k-1}, u_k)$ $\in$
$\{0,1\}^* \times \ldots \times \{0,1\}^*$ is encoded to
 \ ${\sf code}(u_1 \, \# \ \ldots \ u_{k-1} \, \#) \ u_k$ $=$
${\sf code}(u_1) \ 11 \ \ldots \ {\sf code}(u_{k-1}) \ 11 \ u_k$ 
$\ \in \ \{0,1\}^*$. 
A function $f \in {\sf fP}$ is encoded to $f^C \in {\cal RM}_2^{\sf P}$, 
defined by ${\sf domC}(f^C) = {\sf code}({\sf Dom}(f) \, \#)$, so
 \ ${\sf Dom}(f^C) = {\sf code}({\sf Dom}(f)) \ 11 \, \{0, 1\}^*$; and

\smallskip

 \ \ \ \ \ $f^C({\sf code}(x \, \#) \, v) \ = \ {\sf code}(f(x) \ \#) \ v$, 

\smallskip

\noindent for all $x \in {\sf Dom}(f)$ and $v \in \{0,1\}^*$; equivalently,
$f^C({\sf code}(x) \ 11 \ v) = {\sf code}(f(x)) \ 11 \ v$.
Then for every $L \subseteq \{0,1\}^*$, ${\sf code}(L \#)$ is a prefix code, 
which belongs to {\sf P} iff $L$ is in {\sf P}.
And $f \in {\sf fP}$ iff $f^C \in {\cal RM}_2^{\sf P}$.
The transformation $f \mapsto f^C$ is a isomorphic embedding of {\sf fP}
into ${\cal RM}_2^{\sf P}$; moreover, $f^C$ is regular in 
${\cal RM}_2^{\sf P}$ iff $f$ is regular in {\sf fP}. 
From here on, the alphabet denoted by $A$ will always be $\{0,1\}$. 

\smallskip

In \cite{s1f} (Section 4) we introduced a notion of {\em polynomial program}
for Turing machines with built-in polynomial counter (for input balance and
time-complexity). These programs form a machine model that characterizes the 
functions in {\sf fP}. 
For a polynomial program $w$, we let $\phi_w \in {\sf fP}$ denote the 
function computed by this program. 
For every polynomial $q$ of the form $q(n) = a \, n^k + a$ (where $a, k$ are
positive integers), we constructed an {\em evaluation map} 
${\sf ev}_{q}^C \in {\sf fP}$ such that for every polynomial program $w$ 
with built-in polynomial $p_w(n) \le q(n)$ (for all $n \ge 0$), and all 
$x \in A^*$,

\smallskip

 \ \ \ \ \ ${\sf ev}_{q}^C\big({\sf code}(w) \ 11 \ x\big)$ \ $=$
      \ ${\sf code}(w) \ 11 \ \phi_w(x)$

\smallskip

\noindent if $x \in {\sf Dom}(\phi_w)$; if $x \not\in {\sf Dom}(\phi_w)$
then ${\sf ev}_{q}^C\big({\sf code}(w) \, 11 \, x\big)$ is undefined. 
We used ${\sf ev}_{q}^C$, with any polynomial $q$ of degree $\ge 2$ with
large enough coefficient, to prove the following: 
First, {\sf fP} is finitely generated (Theorem 4.5 in \cite{s1f}).  
Second, ${\sf ev}_{q}^C$ is complete in {\sf fP} with respect to inversive 
polynomial reduction (Section 5 of \cite{s1f}). Later in this paper (Def.\
\ref{simulationDEF} and following) we define completeness and various 
reductions for ${\cal RM}_2^{\sf P}$, along the same lines as for {\sf fP}. 

Note that {\sf fP} and ${\cal RM}_2^{\sf P}$, in their entirety, do not have
evaluation maps that belong to {\sf fP}, respectively ${\cal RM}_2^{\sf P}$
(since such maps would not have polynomially bounded complexity). That is
the reason why we restrict {\sf ev} and {\sf evR} to complexity $\le q(.)$,
and why we need precise machine models for {\sf fP} and
${\cal RM}_2^{\sf P}$ (as opposed to more intuitive ``higher-level'' 
models).

\medskip

In Section 2 we define a machine model that characterizes the functions in 
${\cal RM}_2^{\sf P}$; and for any large enough polynomial $q$ we construct 
evaluation maps ${\sf evR}^C_q$ and ${\sf evR}^{CC}_q$ for the functions in 
${\cal RM}_2^{\sf P}$ that have balance and time-complexity $\le q$. 
We prove that ${\sf evR}^{CC}_q$ is complete in ${\cal RM}_2^{\sf P}$ (and 
in {\sf fP}) with respect to inversive Turing reduction.
In Section 3 we prove that ${\cal RM}_2^{\sf P}$ is not finitely generated,
and in Section 4 we show that infinite generation has some complexity 
consequences, i.e., infinite generation can be used for a time-complexity
lower-bound argument.

%%%%%%%%%%%%%%%%%%%%%
%%% Section
%%%%%%%%%%%%%%%%%%%%%
\section{Machine model and evaluation maps for ${\cal RM}_2^{\sf P}$}

The evaluation map $\, {\sf ev}^C_q: {\sf code}(w) \ 11 \ x$ $\longmapsto$ 
${\sf code}(w) \ 11 \ \phi_w(x)$, that we constructed for {\sf fP} in 
\cite{s1f}, works in particular when $\phi_w$ $\in$ ${\cal RM}_2^{\sf P}$ 
(provided that $\phi_w$ has time-complexity and input-balance $\le q$). 
But ${\sf ev}^C_q$ is not a right-ideal morphism and, moreover,
${\sf ev}^C_q$ can evaluate functions that are not in ${\cal RM}_2^{\sf P}$.
We want to construct an evaluation map that belongs to 
${\cal RM}_2^{\sf P}$, and that evaluates exactly the elements of 
${\cal RM}_2^{\sf P}$ that have balance and complexity $\le q$.
In \cite{s1f} we constructed a machine model for {\sf fP}, namely a class of 
Turing machines with built-in polynomial counter (for controlling the 
time-complexity and the input-balance). We will refine these Turing machines 
in order to obtain a machine model for accepting the right ideals in {\sf P}, 
and for computing the functions in ${\cal RM}_2^{\sf P}$.

We will consider deterministic multi-tape Turing machines with input-output
alphabet $A$, with a read-only input tape, and a write-only output tape. 
Moreover we assume that on the input tape and on the output tape, the head 
can only move to the right, or stay in place (but cannot move left). 
We assume that the input tape has a left endmarker $\#$, and a right 
endmarker {\sf B} (the blank symbol).  
At the beginning of a computation of such a machine $M$ 
on input $z \in A^*$, the input tape has content $\# \, z \, {\sf B}$, with 
the input tape head on $\#$ ; initially, all other tapes are blank (i.e., 
they are filled with infinitely many copies of the letter {\sf B}). 
The output tape does not need endmarkers (since it is write-only). 
We assume that $M$ has a special {\em output state} $q_{\rm out}$, and 
that $M$ only goes to state $q_{\rm out}$ when the output is complete; 
the output state is a halting state (i.e., $M$ has no transition from 
state $q_{\rm out}$).
An important convention for a Turing machine $M$ with non-total input-output 
function $f_M$ is the following: 
{\em If $M$ on input $x$ halts in a state that is not $q_{\rm out}$, then 
there is no output (even if the output tape contains a non-blank word).}
So, in that case, $f_M(x)$ is undefined.  
The content of the output tape is considered unreadable, or hidden, until 
the output state $q_{\rm out}$ is reached. 

This kind of Turing machine can compute any partial recursive function (the 
restrictions on the input and output tapes do not limit the machine, because
of the work-tapes).
To compute a function in {\sf fP}, we add a built-in polynomial (used as a 
bound on input balance and time-complexity); see Section 4 in \cite{s1f}.

In order to obtain a machine model for the functions in 
${\cal RM}_2^{\sf P}$ the above Turing machines (with built-in polynomial) 
will be restricted so that they compute right-ideal morphisms of $A^*$.
This is done in two steps: First, sequential functions and sequential Turing 
machines are introduced. From this it is easy to obtain a class of Turing
machines that compute right-ideal morphisms (which are a special kind of 
sequential functions).
Recall that by ``function'' we mean partial function.  
By definition, a function $f: A^* \to A^*$ is {\em sequential} iff 

\smallskip 

 \ \ \ {\em for all $x_1, x_2 \in {\sf Dom}(f)$: \ if $x_1$ is a prefix of 
  $x_2$ then $f(x_1)$ is a prefix of $f(x_2)$. }

\smallskip

\noindent Obviously, every right-ideal morphism is a sequential function.

A {\em sequential Turing machine} is a deterministic multi-tape 
Turing machine $M$ (with special input tape and special output tape and
output state, according to the conventions above), with input-output 
function $f_M$, such that the following holds.

\smallskip

 \ \ \ {\em For every $x \in {\sf Dom}(f_M)$ and every word $z \in A^*$:
  \ in the computation of $M$ on input $xz$, 

 \ \ \ the input-tape head does not start reading 
 \ $z \, {\sf B}$ \ until $f_M(x)$ has been written on the output tape. }

\smallskip

\noindent To ``read a letter $\ell$'' (in $z {\sf B}$) means to make a 
transition whose input letter is $\ell$.  
So, the input tape has content $\# \, xz \, {\sf B}$, with the input-tape 
head on the left-most letter of $z {\sf B}$ (but no transition has been made 
on that letter yet), and the output tape now has content $f_M(x)$. 
Of course, at this moment the computation of $M$ on input $xz$ is not 
necessarily finished; the state is not necessarily $q_{\rm out}$, the 
output might still grow, and $q_{\rm out}$ might be reached 
eventually, or not; if $q_{\rm out}$ is never reached, there is no 
final output.

The sequential Turing machines form a machine model for the partial 
recursive sequential functions. If we let the machines have a built-in 
polynomial we obtain a machine model for the sequential functions in 
{\sf fP}.

Finally, to obtain a machine model for the functions in ${\cal RM}_2^{\sf P}$ 
we take the sequential Turing machines with built-in polynomial, 
with the following additional condition.

\smallskip

 \ \ \ {\em For every $x \in {\sf Dom}(f_M)$ and every word $z \in A^*$:
  \ in the computation of $M$ on input $xz$,
 
 \ \ \ once $f_M(x)$ has been written on the output tape (after $x$ was read
  on the input tape), 

 \ \ \ the remaining input $z$ is copied to the output tape;
 at this point the state $q_{\rm out}$ is reached. } 

\smallskip

\noindent We call such a machine an {\em ${\cal RM}_2^{\sf P}$-machine}. 

\smallskip

The following shows how, from an {\sf fP}-machine for a function $f$, an 
${\cal RM}_2^{\sf P}$-machine for $f$ can be constructed, provided that
$f \in  {\cal RM}_2^{\sf P}$.  

Let us first consider right ideals in {\sf P}, rather than functions. For 
any polynomial program $w$ for a Turing machine $M_w$ that accepts a 
language $L \in {\sf P}$, we construct a new polynomial program $v$ 
describing a Turing machine $M_v$ that behaves as follows: 
On any input $x \in \{0,1\}^*$, $M_v$ successively examines prefixes of $x$ 
until  it finds a prefix, say $p$, 
that is accepted by $M_w$; $M_v$ does not read the letter of $x$ that comes 
after $p$ until it has decided that $p \not\in L$. As soon as 
$M_v$ finds a prefix $p$ of $x$ such that $p \in L$, $M_v$ accepts the whole 
input $x$.
If $M_w$ accepts no prefix of $x$, $M_v$ rejects $x$. Thus, $M_v$ accepts
$LA^*$ (the right ideal generated by $L$); if $L$ is a right ideal then 
$L A^* = L$.
If $M_w$ has time-complexity $\le T(.)$ (a polynomial) then $M_v$ has 
time-complexity $\le T(.)^2$. 
% Binary search for a prefix does not seem to work, unless the given TM 
% already accepts a right ideal.

Let us now consider functions in ${\cal RM}_2^{\sf P}$. 
Given any polynomial program $w$ for a function $\phi_w \in {\sf fP}$, we 
construct a new polynomial program $v$ such that $M_v$, on input $x$, 
successively examines all prefixes of $x$ until it finds a prefix $p$ in 
${\sf Dom}(\phi_w)$; let $\phi_w(p) = y$. 
Then, on input $x$, the machine $M_v$ outputs $y \, z$, where $z$ is such 
that $x = p \, z$.  Note that since $p$ is the shortest prefix of $x$ such 
that $p \in {\sf Dom}(\phi_w)$, we actually have 
$p \in {\sf domC}(\phi_w)$ (if ${\sf Dom}(\phi_w)$ is a right ideal).    
The machine $M_v$ does not read the letter of $x$ that comes after a prefix
$p$ until it has decided that $p \not\in {\sf Dom}(\phi_w)$ or
$p \in {\sf domC}(\phi_w)$.
Hence, the function computed by $M_v$ is in ${\cal RM}_2^{\sf P}$.  
This construction describes a transformation
$f \in {\sf fP} \longmapsto f_{\sf pref} \in {\cal RM}_2^{\sf P}$, where
$f_{\sf pref}$ is defined as follows: 

\smallskip

\hspace{.5in}  $f_{\sf pref}(x) \ = \ f(p) \ z$, 

\smallskip

\noindent where $x = p \, z$, and $p$ is the shortest prefix of $x$ that 
belongs to ${\sf Dom}(f)$; so, $p \in {\sf domC}(f_{\sf pref})$.  Thus for
every $f \in {\sf fP}$ we have: $f \in {\cal RM}_2^{\sf P}$ \, iff \, 
$f_{\sf pref} = f$.  

\smallskip

Based on ${\cal RM}_2^{\sf P}$-machines we can construct evaluation maps 
for ${\cal RM}_2^{\sf P}$.
Let $q$ be a polynomial where $q(n) = a \, n^k + a$ for some integers 
$a, k \ge 1$.  We define ${\sf evR}^C_{q}$, as follows:

\medskip

 \ \ \ \ \  ${\sf evR}_q^C \big({\sf code}(w) \ 11 \ x\big)$ \  $=$
 \ ${\sf code}(w) \ 11 \ \phi_w(x)$,

\medskip

\noindent 
for all ${\cal RM}_2^{\sf P}$-programs $w$ with built-in polynomial
$p_w \le q$, and for all $x \in {\sf Dom}(\phi_w)$. 
The details of the construction are the same as for ${\sf ev}^C_{q}$; see 
Section 4 in \cite{s1f}. 
Although ${\sf evR}^C_{q}$ belongs to ${\cal RM}_2^{\sf P}$ and evaluates 
all ${\cal RM}_2^{\sf P}$-programs $w$ with built-in polynomial $\le q$, 
we will prove in Theorem \ref{RMqFinGen} that the complexity of 
${\sf evR}^C_{q}$ is higher than $q$.

The following {\em doubly coded evaluation function} is usually more useful
for ${\cal RM}_2^{\sf P}$-programs.  It is defined by 

\medskip

 \ \ \ \ \
${\sf evR}_q^{CC} \big( {\sf code}(w) \ 11 \ {\sf code}(u) \ 11 \ v \big)$
 \ $=$ \   ${\sf code}(w) \ 11 \ {\sf code}(\phi_w(u)) \ 11 \ v$,

\medskip

\noindent when $u \in {\sf domC}(\phi_w)$, $v \in A^*$, and $w$ is as before.

\medskip

To give a relation between ${\sf evR}_q^C$ and ${\sf evR}_q^{CC}$ we will 
use the following partial recursive right-ideal morphism $\gamma$,
defined for very ${\cal RM}_2^{\sf P}$-program $w \in A^*$ and every
$x \in {\sf Dom}(\phi_w)$ by

\medskip

 \ \ \ \ \ \ $\gamma\big({\sf code}(w) \ 11 \, x \big) \ = \ $
${\sf code}(w) \ 11 \ {\sf code}(u) \ 11 \, v$,

\medskip

\noindent where $x = u v$, and $u$ is the shortest prefix of $x$ such that
$u \in {\sf Dom}(\phi_w)$;  equivalently, $u \in {\sf domC}(\phi_w)$.
When $x \not \in {\sf Dom}(\phi_w)$, 
$\gamma\big({\sf code}(w) \ 11 \, x \big)$ is undefined. 
Essentially, $\gamma$ finds the shortest prefix of $x$ that belongs to
${\sf Dom}(\phi_w)$ (or equivalently, to ${\sf domC}(\phi_w)$).
The function $\gamma$ can be evaluated by examining successively longer
prefixes of $x$ until a prefix $u \in {\sf Dom}(\phi_w)$ is fund.
So $\gamma$ is computable with recursive domain, when $w$ ranges over 
${\cal RM}_2^{\sf P}$-programs.

For any fixed ${\cal RM}_2^{\sf P}$-program $w$, let $\gamma_w$ be 
$\gamma$ restricted to this $w$, i.e.,
$\gamma_w = \gamma|_{{\sf code}(w) \, 11 \, A^*}$. In other words,
${\sf Dom}(\gamma_w) =  {\sf code}(w) \ 11 \ {\sf Dom}(\phi_w)$, and

\medskip

 \ \ \ \ \ \ $\gamma_w({\sf code}(w) \, 11 \, uv)$ $ \ = \ $ 
${\sf code}(w) \ 11 \ {\sf code}(u) \ 11 \, v$

\medskip

\noindent when $u \in {\sf domC}(\phi_w)$, $v \in A^*$.  
Similarly we define $\gamma^o_w$ by 
${\sf Dom}(\gamma^o_w) =  {\sf Dom}(\phi_w)$ (as opposed to
${\sf code}(w) \ 11 \ {\sf Dom}(\phi_w)$), and

\medskip

 \ \ \ \ \ \ $\gamma^o_w(uv) $ $ \ = \ $  
${\sf code}(w) \ 11 \ {\sf code}(u) \ 11 \, v$

\medskip

\noindent when $u \in {\sf domC}(\phi_w)$, $v \in A^*$. So,
${\sf Im}(\gamma_w^o) = {\sf Im}(\gamma_w)$ $= $ 
${\sf code}(w) \ 11 \ {\sf domC}(\phi_w) \ 11 \, A^*$.

Then $\gamma_w$ and $\gamma^o_w$ belong to ${\cal RM}_2^{\sf P}$ for 
every fixed $w$. But $\gamma$ itself is not polynomial-time computable, since 
it has to work for all possible ${\cal RM}_2^{\sf P}$-programs $w$.

\smallskip

Another restricted form of $\gamma$ that belongs to ${\cal RM}_2^{\sf P}$ is
obtained by choosing a fixed polynomial $q$, and defining $\gamma_q$ as the 
restriction of $\gamma$ to the set

\smallskip

 \ \ \ \ \ \ $\{ {\sf code}(w) \, 11 \, x \ :$
 \ $w$ is a ${\cal RM}_2^{\sf P}$-program with built-in 
               polynomial $\le q$, and $x \in {\sf Dom}(\phi_w) \}$.

\smallskip

\noindent Hence, $\gamma_q \in {\cal RM}_2^{\sf P}$. 

\medskip

We also define the functions $\pi_0$, $\pi_1$, $\rho_0$, $\rho_1$
$\in$ ${\cal RM}_2^{\sf P}$ by 
 \, $\pi_a(x) = ax$, 
 \, $\rho_a(ax) = x$,  
for all $x \in \{0,1\}^*$ and $a \in \{0,1\}$. 
For a word $w = a_n \ldots a_1$ with $a_i \in \{0,1\}$ we denote 
$\, \pi_{a_n} \circ \ \dots \ \circ \pi_{a_1} \,$ by $\pi_w$, and 
$\, \rho_{a_n} \circ \ \dots \ \circ \rho_{a_1} \,$ by $\rho_w$. 

\medskip

Then we have: 
 \ $\gamma^o_w = \gamma_w \circ \pi_{{\sf code}(w) \, 11}$, and 
 \ $\gamma_w = \gamma^o_w \circ \rho_{{\sf code}(w) \, 11}$.

\medskip
 
Another important function in ${\cal RM}_2^{\sf P}$ is the decoding 
function, defined for any $u, v \in A^*$ by 

\medskip

 \ \ \ \ \ \
${\sf decode}({\sf code}(u) \ 11 \, v) \ = \ uv$,

\medskip

\noindent so ${\sf domC}({\sf decode}) = \{00, 01\}^* \, 11$, and
${\sf imC}({\sf decode}) = \{\varepsilon\}$.
We also define a second-coordinate decoding function, for all 
$u_1, u_2, v \in A^*$, by

\medskip

 \ \ \ \ \
${\sf decode}_2\big({\sf code}(u_1) \ 11 \ {\sf code}(u_2) \ 11 \ v \big)$
$\, = \, $  ${\sf code}(u_1) \ 11 \ u_2 \ v$.

\medskip

\noindent So, ${\sf decode}_2 \in {\cal RM}_2^{\sf P}$, 
 \ ${\sf domC}({\sf decode}_2) = \{00,01\}^* \, 11 \, \{00,01\}^* \, 11$,
and ${\sf imC}({\sf decode}_2) = \{00,01\}^* \, 11$.

\medskip

Now we can formulate a relation between ${\sf evR}_q^C$ and 
${\sf evR}_q^{CC}$:

\medskip

 \ \ \ \ \  ${\sf evR}_q^C$ $=$
${\sf decode}_2 \circ {\sf evR}_q^{CC} \circ \gamma_q$.

\medskip

In order to show that ${\sf evR}_q^{CC}$ is complete with respect to 
inversive reduction in ${\cal RM}_2^{\sf P}$, we will adapt the padding and 
unpadding functions (defined for {\sf fP} in \cite{s1f}, Section 4) to 
${\cal RM}_2^{\sf P}$.  Although for ${\cal RM}_2^{\sf P}$ we keep the same 
names as for the corresponding (un)padding functions in {\sf fP}, the 
functions are slightly different.  The padding procedure begins with the 
function ${\sf expand}(.)$, defined by

\medskip

 \ \ \ \ \       
${\sf expand}\big({\sf code}(w) \ 11 \ {\sf code}(u) \ 11 \ v \big)$

\smallskip

 \ \ \ \ \   
$ \ = \ $
${\sf code}({\sf ex}(w)) \ 11 $
  $0^{4\, |{\sf code}(u)|^2 + 8 \, |{\sf code}(u)| + 2} \ 01$
  $ \, {\sf code}(u) \ 11 \ v$, 

\medskip

\noindent for all $u \in {\sf domC}(\phi_w)$, $v \in A^*$, and 
${\cal RM}_2^{\sf P}$-programs $w$.  The word 
$\, 0^{4\, |{\sf code}(u)|^2 + 8 \, |{\sf code}(u)| + 2} \ 01 \,$ is of the 
form ${\sf code}(s)$ for a word $s \in 0^* 1$; and 
$\, 0^{4\, |{\sf code}(u)|^2 + 8\, |{\sf code}(u)| + 2} \ 01$
${\sf code}(u) \, $ is also a code word, namely ${\sf code}(s u)$.   
Since $0^*1$ and its subset $(00)^* \, 01$ are prefix codes, 
${\sf code}(s) =$
$0^{4\, |{\sf code}(u)|^2 + 8\, |{\sf code}(u)| + 2} \ 01 \, $ 
is uniquely determined as a prefix of ${\sf code}(s u)$.

Here, ${\sf ex}(w)$ is an ${\cal RM}_2^{\sf P}$-program obtained from $w$ 
so that

\medskip

 \ \ \ \ \      
$\phi_{{\sf ex}(w)}\big((00)^h \ 01 \ {\sf code}(u) \ 11 \, v \big)$ 
$ \ = \ $
$(00)^h \ 01 \ {\sf code}(\phi_w(u)) \ 11 \ v$,

\medskip

\noindent for all $u \in {\sf domC}(\phi_w)$, $v \in A^*$, and $h > 0$.
Moreover, if $n \mapsto a \, n^k + a$ is the built-in polynomial of the 
program $w$ then the built-in polynomial of ${\sf ex}(w)$ is 

\medskip

 \ \ \ \ \          
$p_e(n) = a_e \, n^{\lceil k/2 \rceil} + a_e$, \ with \
$a_e \, = \, \max\{12, \, \lceil a/2^k \rceil +1\}$.

\medskip

\noindent The detailed justification of the numbers used in the definition 
of {\sf expand} and {\sf ex} (as well as {\sf reexpand}, {\sf recontr}, and 
{\sf contr} below) is given in \cite{s1f}, Section 4.

It is important that ${\sf expand}$ uses the prefix $u$ of $x$ for padding  
(in the format ${\sf code}(u) \, 11$, where $u \in {\sf domC}(\phi_w)$).  
If the whole input $x$ were used for computing the amount of padding, 
${\sf expand}$ would not be a right-ideal morphism. This is the reason why 
we introduce $\gamma_w$ or $\gamma^o_w$, in order to isolate the prefix 
$u \in {\sf domC}(\phi_w)$ of $x$.

We iterate expansion (padding) by applying the following function, where 
${\sf ex}(.)$ is as above:

\medskip

 \ \ \ \ \    
${\sf reexpand}\big({\sf code}({\sf ex}(z)) $
$ \,  11 \ 0^k \ 01 \ {\sf code}(u) \ 11 \ v \big)$

\smallskip

 \ \ \ \ \   
$ \ = \ $   ${\sf code}({\sf ex}(z)) \ 11$
 $0^{4 k^2 + 8 k + 2} \ 01 \ {\sf code}(u) \ 11 \ v$, 

\medskip

\noindent where $k > 0$, $u, v \in A^*$, and $z$ is any 
${\cal RM}_2^{\sf P}$-program; $k$ is even in the context where 
{\sf reexpand} will be used.

Repeated contraction (unpadding) is carried out by applying the following 
function, for $k > 0$: 

\medskip

 \ \ \ \ \    
${\sf recontr}\big({\sf code}({\sf ex}(z)) \ 11 $
 $(00)^k \ 01 \ {\sf code}(y) \ 11 \ v \big)$

\smallskip

 \ \ \ \ \   $ \ = \ $    ${\sf code}({\sf ex}(z)) \ 11 $
 $(00)^{\max\{1, \ \lfloor \sqrt{k}/2 \rfloor -1\}} \ 01$
 ${\sf code}(y) \ 11 \ v$;

\medskip

\noindent note that 
$\, \max\{1, \ \lfloor \sqrt{k}/2 \rfloor -1\} \, \ge \, 1$.

The unpadding procedure ends with the application of the function

\medskip

 \ \ \ \ \     
${\sf contr}\big({\sf code}({\sf ex}(z)) \ 11 $
 $(00)^k \ 01 \ {\sf code}(y) \ 11 \ v \big)$ 
$ \ = \ $
${\sf code}(z) \ 11 \ {\sf code}(y) \ 11 \ v$,

\medskip

\noindent if $\, 2 \, \le \, |(00)^k| \, = \, 2 k \, \le $
$4\, |{\sf code}(y)|^2 + 8\, |{\sf code}(y)| + 2$.

\smallskip

The functions ${\sf expand}(.)$, ${\sf ex}(.)$, ${\sf reexpand}(.)$, 
${\sf recontr}(.)$, and ${\sf contr}(.)$, are undefined in the cases where 
no output has been specified above.

%%%%%%%%%%%
\begin{lem} \label{RMsimulation} \     
Let $q_2$ be the polynomial defined by $q_2(n) = 12 \, n^2 + 12$.
For any $\phi_w \in {\cal RM}_2^{\sf P}$, where $w$ is a 
${\cal RM}_2^{\sf P}$-program with built-in polynomial $q$ (of the form 
$q(n) = a \, n^k + a$ for positive integers $a, k$), we have for all 
$u \in {\sf domC}(\phi_w)$, $v \in A^*$:

\medskip

\noindent $(\star)$
\hspace{.4in}
$\phi_w(u v)$

\smallskip

\hspace{.4in} $ \ = \ $ 
${\sf decode} \circ \rho_{{\sf code}(w) \, 11}$ $\circ$
${\sf contr} \circ {\sf recontr}^{2m} \circ {\sf evR}_{q_2}^{CC}$
$\circ$ ${\sf reexpand}^{m} \circ {\sf expand} \circ \gamma^o_w(u v)$

\smallskip

\hspace{.4in} $ \ = \ $  
$\rho_{{\sf code}(w) \, 11} \circ {\sf decode}_2$ $\circ$
${\sf contr} \circ {\sf recontr}^{2m} \circ {\sf evR}_{q_2}^{CC}$
$\circ$ ${\sf reexpand}^{m} \circ {\sf expand} \circ \gamma^o_w(u v)$, 

\medskip

\noindent where $m = \lceil \log_2 (a + k) \rceil$.
\end{lem}
{\bf Proof.} This is similar to the proof of Prop.\ 4.5 in \cite{s1f}, with
a few modifications. For $u \in {\sf domC}(\phi_w)$, $v \in A^*$,

\medskip

$u v \ \ \stackrel{\gamma^o_w}{\longmapsto} \ \ $
${\sf code}(w) \ 11 \ {\sf code}(u) \ 11 \ v $
 
\medskip

$\stackrel{\sf expand}{\longmapsto} \ \ $
${\sf code}({\sf ex}(w)) \ 11 $
  $0^{4 \, |{\sf code}(u)|^2 + 8\, |{\sf code}(u)| +2} \ 01 $ 
  ${\sf code}(u) \ 11 \ v$

\medskip

$\stackrel{{\sf reexpand}^{m}}{\longmapsto} \ \ $
${\sf code}({\sf ex}(w)) \ 11 \ 0^{N_{m+1}} \ 01 $
      ${\sf code}(u) \ 11 \ v$,

\medskip

\noindent where $N_1 = 4\, |{\sf code}(u)|^2 + 8\, |{\sf code}(u)| +2$, so
$|0^{N_1} \, 01| = (2 \, (|{\sf code}(u)| + 1))^2$;     
by induction,
$N_i = 4\, N_{i-1}^2 + 8\,  N_{i-1} + 2$ \ for $1 < i \leq 2m+1$,
and $|0^{N_i} \, 01| = (2 \, (N_{i-1} + 1))^2$.
The above string, which will now be the argument of ${\sf ev}_{q_2}^{CC}$, 
has length $> N_{m+1} + 2 + |{\sf code}(u)|$, which is much larger than the 
time it takes to simulate the machine with program $w$ on input $u$. So 
${\sf evR}_{q_2}^{CC}$ can now be applied correctly. 
Continuing the calculation,

\medskip

$\stackrel{{\sf evR}_{q_2}^{CC}}{\longmapsto} \ \ $
${\sf code}({\sf ex}(w)) \ 11 \ 0^{N_{m+1}} \ 01 \ {\sf code}(\phi_w(u))$
       $11 \ v$
 
\medskip

$\stackrel{{\sf recontr}^{2m}}{\longmapsto} \ \ $
${\sf code}(w) \ 11 \ 00 \ 01 \ {\sf code}(\phi_w(u))$
       $ 11 \ v $.

\medskip

\noindent We use $2m$ in ${\sf recontr}^{2 m}$ because $\phi_w(u)$ could be 
much shorter than $u$; but because of polynomial input balance, 
$|u| \le p_w(|\phi_w(u)|)$.
Note that doing more input padding than necessary does not do any harm;
and recontracting (unpadding) more than needed has no effect (by the 
definition of ${\sf recontr}$).
Hence {\sf contr} can now be applied correctly. We complete the calculation: 

\smallskip

$\stackrel{\sf contr}{\longmapsto} \ \ $
${\sf code}(w) \ 11 \ {\sf code}(\phi_w(u)) \ 11 \ v$
 \ \ $\stackrel{{\sf decode}_2}{\longmapsto} \ \ $
${\sf code}(w) \ 11 \ \phi_w(u) \ v$
 \ \ $\stackrel{\rho_{{\sf code}(w) \, 11}}{\longmapsto} \ \ $
$\phi_w(u) \ v$.
 \ \ \ $\Box$

\begin{lem} \label{infGenSets} \hspace{-.13in} {\bf .} 
 \ ${\cal RM}_2^{\sf P}$ has the following {\em infinite generating set}:

\smallskip

 \ \ \          
$\{{\sf decode}, \, \rho_0, \, \rho_1, \, \pi_0, \, \pi_1, \, {\sf contr},$
$ \, {\sf recontr}, \, {\sf evR}_{q_2}^{CC}, \, $ 
${\sf reexpand}, \, {\sf expand}\}$

\smallskip

 \ \ \ \ $\cup \ \{ \gamma_w: \, w$ is any ${\cal RM}_2^{\sf P}$-program$\}$.

\smallskip

\noindent Here, ${\sf decode}$ can be replaced by ${\sf decode}_2$.
Yet another infinite generating set of ${\cal RM}_2^{\sf P}$ is \,

\smallskip

 \ \ \   
$\{\rho_0, \rho_1, \pi_0, \pi_1\}$ $ \ \cup \ $
$\{ {\sf evR}_q^C : \, q$ is any polynomial of the form
$q(n) = a \, n^k + a$ with $a, k \in {\mathbb N}_{\ge1}$\}.
\end{lem}
{\bf Proof.} The first infinite generating set follows from Lemma 
\ref{RMsimulation}.  
Recall that $\gamma^o_w = \gamma_w \circ \pi_{{\sf code}(w) \, 11}$.
The second generating set follows in a straightforward way from the proof 
of Prop.\ 4.5 in \cite{s1f}.  
 \ \ \ $\Box$

\begin{pro} \label{regGenSet}
 \ ${\cal RM}_2^{\sf P}$ is generated by a set of {\em regular} elements 
of ${\cal RM}_2^{\sf P}$.
\end{pro}
{\bf Proof.} The generators $\rho_0, \rho_1, \pi_0, \pi_1$ are easily 
seen to be regular. Thus, using the second infinite generating set in Lemma 
\ref{infGenSets}, it is enough to factor ${\sf evR}_q^C$ into regular 
elements. We have: 

\smallskip

 \ \ \ \ \ ${\sf evR}_{q}^C = \rho_{2, q} \circ E_{q}$, 

\smallskip

\noindent where $E_{q}$ and $\rho_{2, q}$ are defined as follows:
For every ${\cal RM}_2^{\sf P}$-program $w$ with built-in polynomial 
$\leq q$, and every $u \in {\sf domC}(\phi_w)$ and $v \in A^*$,

\smallskip

 \ \ \ \ \ $E_{q}\big({\sf code}(w) \ 11 \ {\sf code}(u) \ 11 \ v \big)$
$ \ = \ $
${\sf code}(w) \ 11 \ {\sf code}(u) \ 11 \ {\sf code}(\phi_w(u)) \ 11 \ v$; 

\smallskip

\noindent and for all $z, y, x, v \in A^*$ such that $|x| \le q(|y|)$,

\smallskip

 \ \ \ \ \   
$\rho_{2, q} \big({\sf code}(z) \ 11 \ {\sf code}(x) \ 11 \ {\sf code}(y)$
$11 \ v\big)$ $ \ = \ $
${\sf code}(z) \ 11 \ {\sf code}(y) \ 11 \ v$. 

\smallskip

\noindent The functions are undefined otherwise.
It is easy to see that $E_{q}$ and $\rho_{2, q}$ have
polynomial-time inversion algorithms (i.e., they are regular), and belong 
to ${\cal RM}_2^{\sf P}$. 
 \ \ \ $\Box$

\bigskip

We will show now that ${\sf evR}_{q_2}^{CC}$ is complete in
${\cal RM}_2^{\sf P}$ and in {\sf fP}, with respect to a certain ``inversive
reduction''. 
We need to recall some definitions from \cite{s1f} concerning reductions
between functions in {\sf fP} or ${\cal RM}_2^{\sf P}$, and in particular,
reductions that ``preserve one-wayness'' (inversive reductions).

\begin{defn} \label{simulationDEF} 
 \ Let $f_1, f_2: A^* \to A^*$ be two polynomially balanced right-ideal 
morphisms.

\smallskip

\noindent (1) 
We say that $f_2$ {\em simulates} $f_1$ (denoted by $f_1 \preccurlyeq f_2$) 
iff there exist $\alpha, \beta \in {\cal RM}_2^{\sf P}$ such that 
$f_1 = \beta \circ f_2 \circ \alpha$.

\smallskip

\noindent (2)
We have a {\em polynomial-time Turing simulation} of $f_1$ by $f_2$ (denoted 
by $f_1 \preccurlyeq_{\sf T} f_2$) iff $f_1$ can be computed by an oracle
${\cal RM}_2^{\sf P}$-machine that can make oracle calls to $f_2$; such 
oracle calls can, in particular, be calls on the membership problem of 
${\sf Dom}(f_2)$.
%
%\smallskip
%
%\noindent (3) Consider a set $S \subseteq A^*$. By 
%$({\cal RM}_2^{\sf P})^{S}$ we denote the set of polynomially balanced 
%right-ideal morphisms that are computed by oracle 
%${\cal RM}_2^{\sf P}$-machines with oracle set $S$.
%
%\smallskip
%
%\noindent (4)
%We have a {\em weak Turing simulation} of $f_1$ by $f_2$ (denoted by
%$f_1 \preccurlyeq_{\sf wT} f_2$) iff there exist $\beta, \alpha$ such that 
%$f_1 = \beta \circ f_2 \circ \alpha$, where $\beta \in {\cal RM}_2^{\sf P}$
%and $\alpha \in ({\cal RM}_2^{\sf P})^{{\sf Dom}(f_2)}$.
\end{defn}
In the above definition, $f_1, f_2$ need not be polynomial-time computable.

Since ${\cal RM}_2^{\sf P}$ is ${\cal J}^0$-simple (Prop.\ 2.7 in 
\cite{s1f}), every $f_1 \in {\cal RM}_2^{\sf P}$ is simulated by every 
$f_2 \in {\cal RM}_2^{\sf P} - \{0\}$ (for each of the above simulations).

\begin{defn} \label{invSimulationDEF}  {\bf (Inversive reduction).} 
If $\preccurlyeq_{\sf X}$ is a simulation between right-ideal morphisms 
(e.g., as in the previous definition) then the corresponding inversive 
reduction is defined as follows.  We say that $f_1$ {\em inversively 
$X$-reduces} to $f_2$ (denoted by $f_1 \leqslant_{\sf inv,X} f_2$) iff \\  
{\rm (1)} \ $f_1 \preccurlyeq_{\sf X} f_2$, and \\  
{\rm (2)} \ for every inverse $f_2'$ of $f_2$ there exists an inverse $f_1'$
of $f_1$ such that $f_1' \preccurlyeq_{\sf X} f_2'$; here, $f_2'$ and 
$f_1'$ range over all polynomially balanced right-ideal morphisms
$A^* \to A^*$. 
\end{defn}
Note that ${\cal J}^0$-simplicity (Prop.\ 2.7 in \cite{s1f}) does not apply 
for inversive reduction since $f_2', f_1'$ do not range over just 
${\cal RM}_2^{\sf P}$.  
One easily proves the following about polynomially balanced right-ideal
morphisms $f_1, f_2$ (see \cite{s1f}, Section 5):   

{\em
If $f_1 \leqslant_{\sf inv,T} f_2$ and  $f_2 \in {\cal RM}_2^{\sf P}$, 
then $f_1 \in {\cal RM}_2^{\sf P}$;
 \ if, in addition, $f_2$ is regular, then $f_1$ is regular (equivalently,
if, in addition, $f_1$ is one-way, then $f_2$ is one-way).
}

\begin{defn} \label{completenessDEF} \   
A polynomially balanced right-ideal morphism $f_0$ is {\em complete} in a 
set $S$ (of right-ideal morphisms) with respect to an (inversive) reduction 
$\leqslant_{\sf inv,X}$ \ iff \ $f_0 \in S$, and for all $\phi \in S$: 
$\phi \leqslant_{\sf inv,X} f_0$.
\end{defn}

See Section 5 of \cite{s1f} for more details and properties of these 
simulations and reductions; in \cite{s1f} the focus was on {\sf fP},
whereas here we concentrate on ${\cal RM}_2^{\sf P}$. 
The simulations in Def.\ \ref{simulationDEF} are similar to the standard notions of reductions between decision problems.
The concept of inversive reduction was first introduced in \cite{s1f}; it is 
the appropriate notion of reduction between functions when one-wayness is 
to be preserved under upward reduction (and regularity is to be
preserved under downward reduction).  

In the above definitions we only refer to polynomially balanced inverses; 
this is justified by the following Proposition, according to which 
``balanced functions have balanced inverses''.

\begin{pro} \label{polybalInverse}
 \ Suppose $f$ is a right-ideal morphism with balance $\le q(.)$ (where 
$q(.)$ is a polynomial), and $f$ has an inverse $f'_1$ with time-complexity 
$\le T(.)$.  Then $f$ has an inverse $f'$ with balance $\le q$ and 
time-complexity $\, \le T(.) +  c \, q(.)$ (for some constant $c>1$). 
The inverse $f'$ can be chosen as a restriction of $f'_1$.
\end{pro}
{\bf Proof.} Let $f'$ be the restriction of $f'_1$ to the set

\smallskip

\hspace{0.8in}  $\{ y \in {\sf Dom}(f'_1) :$
  $ \ |y| \le q(|f'_1(y)|) \ \ {\rm and} \ \ |f'_1(y)| \le q(|y|)\}$.

\smallskip

\noindent Then $f'$ obviously has balance $\le q$. Note that since $f'_1$ is
an inverse of $f$ we have ${\sf Im}(f) \subseteq {\sf Dom}(f'_1)$.
To show that $f'$ is an inverse of $f$ it is sufficient to check that the
domain of $f'$ contains ${\sf Im}(f)$. 
Let $y = f(x) \in {\sf Im}(f)$ for some $x \in {\sf Dom}(f)$. Then
$f(f'_1(y)) = y$, since $f'_1$ is an inverse.

Checking $|y| \le q(|f'_1(y)|)$: \ $|y| = |f(f'_1(y))| \le q(|f'_1(y)|$;
the inequality holds since $q$ is a balance for $f$ on input $f'_1(y)$.

Checking $|f'_1(y)| \le q(|y|)$: \ $|f'_1(y)| \le q(|f(f'_1(y))|)$ since 
$q$ is a balance for $f$ on input $f'_1(y)$; and 
$q(|f(f'_1(y))|) = q(|y|)$ since $f(f'_1(y)) = y$.

To find a time-complexity bound for $f'$, we first compute $f'_1(y)$ in time 
$\le T(|y|)$; thereby we also verify that $y \in {\sf Dom}(f'_1)$. 
To check whether $y$ is in the domain of $f'$ we first compare $|y|$ and 
$|f'_1(y)|$ in time $\le |y| + 1$.

Checking $|y| \le q(|f'_1(y)|)$: 
If $|y| \le |f'_1(y)|$ then we automatically have $|y| \le q(|f'_1(y)|)$.
If $|y| \ge |f'_1(y)|$ we compute $q(|f'_1(y)|)$ in time 
$O(q(|f'_1(y)|))$ ($\le O(q(|y|))$), by writing the number $|f'_1(y)|$ in 
binary, and then evaluating $q$ (see Section 4 of \cite{s1f} for a similar 
computation). Then we check $|y| \le q(|f'_1(y)|)$ in time $\le |y| + 1$.
Checking $|f'_1(y)| \le q(|y|)$ is done in a similar way, in time 
$\le O(q(|y|)) + |y| + 1$.
 \ \ \ $\Box$

\begin{thm} \label{evR_complete} \  
The map ${\sf evR}_{q_2}^{CC}$ is complete for ${\cal RM}_2^{\sf P}$ with 
respect to inversive Turing reduction.
\end{thm}
{\bf Proof.} Lemma \ref{RMsimulation} provides the following simulation 
of $\phi_w$ by ${\sf evR}_{q_2}^{CC}$:

\medskip

 \ \ \ \ \ $\phi_w \ = \ {\sf decode} \circ \rho_{{\sf code}(w') \, 11}$
$\circ$ ${\sf contr} \circ {\sf recontr}^{2m} \circ {\sf evR}_{q_2}^{CC}$
$\circ$ ${\sf reexpand}^m \circ {\sf expand} \circ \gamma^o_w$ . 

\medskip

\noindent To obtain an inversive Turing simulation, let ${\sf e}'$ be any 
inverse of ${\sf evR}_{q_2}^{CC}$.  Slightly modifying the proof of 
Prop.\ 5.6 in \cite{s1f}, we apply ${\sf e}'$ to any string of the form  
 
\medskip

 \ \ \ \ \  ${\sf code}({\sf ex}(w)) \ 11 \ 0^{N_{m+1}} \ 11 \, $
${\sf code}(p) \ 11 \ z$, 
%%% $p \in \phi_w({\sf domC}(\phi_w))$, 
%%% not necessarily $\in {\sf imC}(\phi_w)$

\medskip

\noindent where $p \in \phi_w({\sf domC}(\phi_w))$, and $z \in A^*$; 
 \ then for any $p \in {\sf imC}(\phi_w)$ ($\subseteq$
$\phi_w({\sf domC}(\phi_w))$), and $z \in A^*$:

\medskip

 \ \ \ \ \    
${\sf e}'\big({\sf code}({\sf ex}(w)) \ 11 \ 0^{N_{m+1}} \ 11$
    ${\sf code}(p) \ 11 \ z \big)$

\smallskip

 \ \ \ \ \    $ \ = \ $
${\sf code}({\sf ex}(w)) \ 11 \ 0^{N_{m+1}} \ 11 \ {\sf code}(t)$
 $ \, 11 \ z$,

\medskip

\noindent for some $t \in \phi_w^{-1}(p) \subseteq {\sf Dom}(\phi_w)$.
Based on ${\sf e}'$ we now construct an inverse $\phi'_w$ of $\phi_w$ such 
that $\phi'_w \preccurlyeq_{\sf T} {\sf e}'$; 
for any $y \in {\sf Im}(\phi_w)$ we define

\medskip

 \ \ \ \ \ $\phi'_w(y) \ = \ $
${\sf decode} \circ \rho_{{\sf code}(w') \, 11}$ $\circ$ 
${\sf contr} \circ {\sf recontr}^{2m}$ $\circ$ ${\sf e}'$ $\circ$
${\sf reexpand}^m \circ {\sf expand} \circ \delta^o_w(y)$.

\medskip

\noindent Here, $\delta^o_w(y)$ is defined by

\medskip

 \ \ \ \ \ $\delta^o_w(y) \ = \ $
${\sf code}(w) \ 11 \ {\sf code}(p) \ 11 \, z$,

\medskip

\noindent when $y = pz$ with $p \in {\sf imC}(\phi_w)$, $z \in A^*$.
So, $\delta^o_w(.)$ is similar to $\gamma^o_w(.)$, except that 
$\delta^o_w(.)$ uses ${\sf imC}(\phi_w)$, whereas $\gamma^o_w(.)$ uses  
${\sf domC}(\phi_w)$. We saw that $\gamma^o_w \in {\cal RM}_2^{\sf P}$; 
but unless ${\sf P} = {\sf NP}$, $\delta^o_w$ will not be in 
${\cal RM}_2^{\sf P}$ in general.

The value $\delta^o_w(y)$ can be computed by an 
${\cal RM}_2^{\sf P}$-machine $M$ that makes oracle calls to 
${\sf Dom}({\sf e}')$ and to ${\sf e}'$ as follows.  
On input $y$, $M$ considers all prefixes of $y$ of increasing lengths, 
$p_1, \ldots, p_k$, until $p_j \in {\sf Im}(\phi_w)$ is found. 
Since $p_j$ is the first prefix in ${\sf Im}(\phi_w)$, we have 
$p_j \in {\sf imC}(\phi_w)$ and 
$\, \delta^o_w(y) = {\sf code}(w) \ 11 \ {\sf code}(p_j) \ 11 \, z$. 
To test for each $p_i$ whether $p_i \in {\sf Im}(\phi_w)$, $M$ pads $p_i$ 
to produce  $\, 0^{N_{m+1}} \, 11 \ {\sf code}(p_i)$;  if 
$p_i \in {\sf Im}(\phi_w)$ then 
$\, {\sf e}'\big({\sf code}({\sf ex}(w)) \ 11 \ \bullet \big) \,$ is 
defined on input $\, 0^{N_{m+1}} \ 11 \ {\sf code}(p_i)$. Thus, if 
$\, {\sf code}\big({\sf ex}(w))\ 11\ 0^{N_{m+1}}\ 11$
 $ \, {\sf code}(p_i)\big)$ $\not\in$ ${\sf Dom}({\sf e}')$, then 
$p_i \not\in {\sf Im}(\phi_w)$.
On the other hand, if 
${\sf code}\big({\sf ex}(w))\ 11\ 0^{N_{m+1}}\ 11$
 $ \, {\sf code}(p_i)\big)$ $\in$ ${\sf Dom}({\sf e}')$, then let 
$t_i \in \phi_w^{-1}(p_i)$ be such that 

\smallskip

 \ \ \ ${\sf e}'\big({\sf code}({\sf ex}(w)) \ 11 \, 0^{N_{m+1}}$
          $11 \ {\sf code}(p_i)\big)$       $\, = \,$
${\sf code}({\sf ex}(w)) \ 11 \ 0^{N_{m+1}} \ 11 \ {\sf code}(t_i)$.

\smallskip

\noindent One oracle call to ${\sf e}'$ yields this, and hence $t_i$. 
Then we can use $\phi_w$ to check whether $t_i \in {\sf Dom}(\phi_w)$;
and this holds iff $p_i \in {\sf Im}(\phi_w)$. 
This way, $M$ can check whether $p_i \in {\sf Im}(\phi_w)$. 
Thus, if $y \in {\sf Im}(\phi_w)$, $M$ will find
$p_j \in {\sf Im}(\phi_w)$.  
When $y \not\in {\sf Im}(\phi_w)$, $M$ produces no output; this doesn't
matter since we do not care how $\phi'_w$ is defined outside of 
${\sf Im}(\phi_w)$.

Once $\delta^o_w(y)$ is known, the remaining simulation

\smallskip

 \ \ \ \ \
${\sf decode} \circ \rho_{{\sf code}(w') \, 11}$ $\circ$ 
${\sf contr} \circ {\sf recontr}^{2m} \circ {\sf e}'$ $\circ$
${\sf reexpand}^m \circ {\sf expand}$ 

\smallskip

\noindent of ${\sf e}'$, applied to 
$\, \delta^o_w(y) = {\sf code}(w) \, 11 \, {\sf code}(p) \, 11 \, z$, 
yields $\phi'_w(y)$.

The function $\phi'_w$ is an inverse of $\phi_w$: Indeed, for 
$x \in {\sf Dom}(\phi_w)$, we have $\phi_w(x) = p z$ for some 
$p \in {\sf imC}(\phi_w)$, $z \in A^*$. Then

\medskip

${\sf reexpand}^m \circ {\sf expand} \circ \delta^o_w(pz) \ = \ $
${\sf code}({\sf ex}(w)) \ 11 \ 0^{N_{m+1}} \ 11 \ {\sf code}(p) \,$
$ 11 \ z$; 

\medskip 

\noindent and applying ${\sf e}'$ then yields

\medskip

${\sf code}({\sf ex}(w)) \ 11 \ 0^{N_{m+1}} \ 11 \ {\sf code}(t)$
 $ \, 11 \ z$,

\medskip

\noindent for some $t \in \phi_w^{-1}(p)$.  Applying 

\medskip

${\sf decode} \circ \rho_{{\sf code}(w') \, 11}$
$\circ$ ${\sf contr} \circ {\sf recontr}^{2m}$ 

\medskip

\noindent now yields $t z$.  Finally, $\phi_w(tz) = p z$, since 
$t \in \phi_w^{-1}(p)$.  So, $\phi_w \phi'_w \phi_w(x)$ $=$
$\phi_w  \phi'_w(p z) = \phi_w(t z) = p z = \phi_w(x)$.
 \ \ \ $\Box$

\bigskip

\noindent We show next that ${\sf evR}_{q_2}^{CC}$ is not only complete for
${\cal RM}_2^{\sf P}$, but for all of {\sf fP}.

\begin{pro}
 \ The map ${\sf evR}_{q_2}^{CC}$ ($\in {\cal RM}_2^{\sf P}$) is complete 
for {\sf fP} with respect to $\leqslant_{\sf inv, T}$.
\end{pro}
{\bf Proof.} By Prop.\ 5.6 in \cite{s1f}, ${\sf ev}_{q_2}^C$ is complete in 
{\sf fP} for inversive simulation. By Prop.\ 5.17 in \cite{s1f},
${\sf ev}_{q_2}^C \leqslant_{\sf inv} ({\sf ev}_{q_2}^C)^C$. 
Moreover, $({\sf ev}_{q_2}^C)^C$ $\leqslant_{\sf inv, T}$ 
${\sf evR}_{q_2}^{CC}$; indeed,
$({\sf ev}_{q_2}^C)^C \in {\cal RM}_2^{\sf P}$ (since $f \mapsto f^C$ maps
into ${\cal RM}_2^{\sf P}$), and we just saw that ${\sf evR}_{q_2}^{CC}$ is 
complete in ${\cal RM}_2^{\sf P}$. Hence ${\sf ev}_{q_2}^C$ 
$\leqslant_{\sf inv}$ $({\sf ev}_{q_2}^C)^C$ $\leqslant_{\sf inv, T}$
${\sf evR}_{q_2}^{CC}$. 
 \ \ \ $\Box$

%%%%%%%%%%%%%%%%%%%%%%%%%%%%%%%%%%%%%%%%%%%%%%%%%%%%%%%%%%
%%% Section
%%%%%%%%%%%%%%%%%%%%%%%%%%%%%%%%%%%%%%%%%%%%%%%%%%%%%%%%%%
\section{Non-finite generation }

In \cite{s1f} we proved that {\sf fP} is finitely generated, and we left 
open the question whether ${\cal RM}_2^{\sf P}$ is also finitely generated.
We will now answer this question negatively.
We will use the following general compactness property: If a semigroup $S$ 
is finitely generated, and if $\Gamma$ is any infinite generating set of $S$, 
then $S$ is generated by some finite subset of this set $\Gamma$.

\begin{thm} \hspace{-.13in} {\bf .} \label{RMnotFinGen} 
%{\bf (Non-finite generation).} 
 \ \ ${\cal RM}_2^{\sf P}$ is not finitely generated.
\end{thm}
{\bf Proof.} We saw that ${\cal RM}_2^{\sf P}$ is generated by the infinite
set

\smallskip

\hspace{0.2in} 
$\{\rho_0, \, \rho_1, \, \pi_0, \, \pi_1, \, {\sf decode}_2, $
$ \, {\sf contr}, \,  {\sf recontr}, \, $
${\sf evR}_{q_2}^{CC}, \, {\sf reexpand}, \, {\sf expand}\}$  

\smallskip

\hspace{0.2in} $\cup \ $
$\{\gamma_w: \, w$ is an ${\cal RM}_2^{\sf P}$-program\}.

\smallskip

\noindent Let us assume, by contradiction, that ${\cal RM}_2^{\sf P}$ is
finitely generated. Then a finite generating set can be extracted from this
infinite generating set, so ${\cal RM}_2^{\sf P}$ is generated by

\smallskip

\hspace{0.2in} $\Gamma_{\sf fin} \ = \ $
$\{\rho_0, \, \rho_1, \, \pi_0, \, \pi_1, \, {\sf decode}_2, \, {\sf contr},$
$ \,  {\sf recontr}, \, {\sf evR}_{q_2}^{CC}, \, {\sf reexpand}, $
$ \, {\sf expand} \}$ $ \ \cup \ $ $\{\gamma_i: i \in F\}$,

\smallskip

\noindent where $F$ is some finite set of ${\cal RM}_2^{\sf P}$-programs.
So for every $\gamma_w$ there is a word in $\Gamma_{\sf fin}^*$ that 
expresses $\gamma_w$ as a finite sequence of generators. Recall that 
${\sf Dom}(\gamma_w)$ $=$ ${\sf code}(w) \ 11 \ {\sf Dom}(\phi_w)$,
and for any $x \in {\sf Dom}(\phi_w)$, 

\smallskip

\hspace{0.2in} $\gamma_w({\sf code}(w) \ 11 \, x)$ \ $=$ \
 \ ${\sf code}(w) \ 11 \ {\sf code}(u) \ 11 \, v$,

\smallskip

\noindent where $x = u v$ and $u \in {\sf domC}(\phi_w)$. 

\smallskip

The proof strategy will consist in showing that there are infinitely many 
functions $\gamma_w$ that do not have a correct representation over 
$\Gamma_{\sf fin}$. 
More precisely, for all ${\cal RM}_2^{\sf P}$-programs $w$ and all 
$u \in {\sf domC}(\phi_w)$, we have $\, \gamma_w({\sf code}(w) \, 11 \, u)$ 
$=$ ${\sf code}(w) \, 11 \, {\sf code}(u) \, 11$; so 
$\, \gamma_w({\sf code}(w) \, 11 \, u)$
$\in$ $\{00,01\}^* \, 11 \, \{00,01\}^* \, 11$. 
On the other hand, we will show that there exist (infinitely many) 
${\cal RM}_2^{\sf P}$-programs $w$ such that for every 
$X \in (\Gamma_{\sf fin})^*$ that represents $\gamma_w$,
there exist (infinitely many) $u \in {\sf domC}(\phi_w)$ such that:
$X({\sf code}(w) \, 11 \, u)$ $=$ 
${\sf code}(w) \, 11 \, {\sf code}(u_1) \, 11 \, u_2$, where $u_2$ is 
non-empty; so, $X({\sf code}(w) \, 11 \, u)$ $\not\in$ 
$\{00,01\}^* \, 11 \, \{00,01\}^* \, 11$. Thus we obtain a contradiction.

\smallskip

We consider the ${\cal RM}_2^{\sf P}$-programs $w$ such that 
${\sf domC}(\phi_w)$ satisfies:

\noindent
(1) \ no word in ${\sf domC}(\phi_w)$ contains 11 as a subsegment;

\noindent
(2) \ for all $i \in F$, \, ${\sf domC}(\phi_i) \neq {\sf domC}(\phi_w)$;

\noindent
(3) \ for any $u \in {\sf domC}(\phi_w)$ and any integer $n>0$, there
exists $v \in {\sf domC}(\phi_w)$ of length $|v| > n$ such that 
$u = u_0 c$, $u_0$ is a prefix of $v$, and $|c| \le 4$.
Equivalently: 

\noindent
$\big(\forall u \in {\sf domC}(\phi_w)\big)\big(\forall n$
$>0 \big) \big(\exists v \in {\sf domC}(\phi_w), $
$ |v| > n\big)\big(\exists u_0, c, z \in A^*\big)$
$[ \, v = u_0 z, \ u = u_0 c, \ |c| \le 4 \, ]$.  
 
\noindent We can picture this as a path in the tree of $A^*$, labeled by 
$u$ and ending at vertex $u$; at vertex $u_0$ along this path, at distance 
$\le 4$ from vertex $u$, a second path branches off and ends at vertex $v$ 
(of length $|v| > n$).

The following family of examples shows that there exist infinitely many 
${\cal RM}_2^{\sf P}$-programs $w$ that satisfy properties (1)-(3). 
In each of these examples (parameterized by $a \in \{0,1\}^*$) we have

\medskip

\hspace{.4in}
${\sf domC}(\phi_w) \ = \ \{{\sf code}(a^n) \ 0010 \ : \ n > 0\}$,

\medskip

\noindent where $a \in \{0,1\}^*$ is any fixed non-empty word (depending on 
$w$), chosen so that ${\sf domC}(\phi_w) \neq {\sf domC}(\phi_i)$ for all 
$i \in F$ (thus property (2) holds). Any word $a$ that is long enough will 
work; indeed, for different words $a$ the above prefix codes are different, 
whereas $F$ is finite.
Property (1) follows from the definition of {\sf code} (namely, 
${\sf code}(0) = 00, \ {\sf code}(1) = 01$). Property (3) holds because
for every $u = {\sf code}(a^m) \ 0010$ and every $n>0$, we can take
$u_0 = {\sf code}(a^m)$ and $v = {\sf code}(a^{n+m}) \, 0010$.
The set $\{{\sf code}(a^n) \ 0010 \ : \ n > 0\}$ is a regular language, with 
regular expression $\, ({\sf code}(a))^+ \, 0010$.
 
\smallskip

Let $X \in \Gamma_{\sf fin}^*$ be a representation of $\gamma_w$, where $w$ 
is any ${\cal RM}_2^{\sf P}$-program from the family of examples above with 
properties (1)-(3).
We will consider certain suffixes $S_i$ of $X$, over $\Gamma_{\sf fin}$.

Let $S_0$ be the shortest suffix of $X$ such that for all 
$u \in {\sf domC}(\phi_w)$, $S_0({\sf code}(w) \, 11 \, u)$ is of
the form $\, {\sf code}(x) \, 11 \, {\sf code}(y) \, 11 \, z$ $\, \in \,$
$\{00,01\}^* \, 11 \, \{00,01\}^* \, 11 \, \{0,1\}^*$.
Then $S_0$ exists since $X$ itself (representing $\gamma_w$) maps 
 \ ${\sf code}(w) \, 11 \, u$ \ to 
 \ ${\sf code}(w) \, 11 \, {\sf code}(u) \, 11 \, \in $
$\, \{00,01\}^* \, 11 \, \{00,01\}^* \, 11 \, \{0,1\}^*$.

Inductively we define $S_1, S_2, \ldots, S_i, \ldots \, $, where $S_i$ is 
the shortest suffix of $X$ that has $S_{i-1}$ as a strict suffix, and 
such that for all $u \in {\sf domC}(\phi_w)$ we have: 

\smallskip

\hspace{.4in} 
$S_i\big({\sf code}(w) \ 11 \, u \big)$ $ \ \in \ $
 $\{00,01\}^* \, 11 \, \{00,01\}^* \, 11 \, \{0,1\}^*$.

\smallskip

\noindent So, $S_i\big({\sf code}(w) \ 11 \, u \big)$ is of the form
$\, {\sf code}(w_1) \, 11 \, {\sf code}(u_1) \, 11 \, u_2 \,$ with 
$w_1, u_1, u_2 \in A^*$.
Then $X = S_N$ for some $N \ge 0$ (and $|X| > N$). 

\smallskip

Theorem \ref{RMnotFinGen} now follows from the next Lemma, according to 
which there are (infinitely many) $u \in {\sf domC}(\phi_w)$ such that
$S_N({\sf code}(w) \ 11 \, u)$ $=$
${\sf code}(w) \ 11 \ {\sf code}(u_1) \ 11 \, u_2$, with $u_2$ 
{\em non-empty}. 
On the other hand, $X$ = $S_N$, and $X$ represents $\gamma_w$, 
hence by the definition of $\gamma_w$ we have for {\em every} 
$u \in {\sf domC}(\phi_w)$:
$S_N({\sf code}(w) \ 11 \, u)$ \, $=$ \, 
${\sf code}(w) \ 11 \ {\sf code}(u) \ 11$; so $u_2$ is empty. Thus, the 
assumption that $X$ (over the finite generating set $\Gamma_{\sf fin}$) 
represents $\gamma_w$, leads to a contradiction.

%%%Lemma
\begin{lem} \label{lemNotFinGen} 
 \ Let $\gamma_w$ be such that 
${\sf domC}(\phi_w) = \{{\sf code}(a^n) \, 0010 : n > 0\}$ for some word 
$a \in \{0,1\}^*$, chosen so that the program $w$ satisfies properties
(1)-(3).  Let $X$ be a word over $\Gamma_{\sf fin}$ that represents 
$\gamma_w$, and let $|X|$ be the length of $X$ over $\Gamma_{\sf fin}$.
Let $S_0, \ldots, S_N$ be the suffixes of $X$ defined above, with $S_N = X$.
  Then there exist $\ell$ and $n$ with $\ell > n > 0$ such that for all
$i = 0, \ldots, N$ and all $u \in {\sf domC}(\phi_w)$ with $|u| \ge \ell$:

\smallskip

 \ \ \ \ \ $S_i({\sf code}(w) \ 11 \, u)$  $=$ 
${\sf code}(w_1) \ 11 \ {\sf code}(u_1) \ 11 \, u_2$, 

\smallskip

\noindent for some $w_1, u_1, u_2$ $\in A^*$. Moreover, $u_2$ has a 
{\em non-empty} common suffix with $u$, and this common suffix has length 
at least $n$.
\end{lem} 
{\bf Proof.} We have for all $u \in {\sf domC}(\phi_w)$:  
$S_i({\sf code}(w) \, 11 \, u)$ $=$
${\sf code}(w_1) \, 11 \, {\sf code}(u_1) \, 11 \, u_2$, for some 
$w_1, u_1, u_2 \in A^*$.  We want to show that there is $\ell$ such that 
for all $u \in {\sf domC}(\phi_w)$ with $|u| \ge \ell$: $u_2$ has a 
non-empty (sufficiently long) suffix in common with $u$; the number $n$ is
an auxiliary parameter.  We take $u$ of the form 
$u = {\sf code}(a^m) \, 0010$ and use induction on $i = 0, \ldots, N$. 

\smallskip

Proof for $S_0$: The only generators from $\Gamma_{\sf fin}$ that can occur 
in $S_0$ are $\pi_0, \pi_1, \rho_0, \rho_1$ and $\gamma_j$ (for $j \in F$). 
Indeed, the other generators in $\Gamma_{\sf fin}$ (namely 
${\sf decode}_2,$ ${\sf contr},$ ${\sf recontr},$ ${\sf evR}_{q_2}^{CC},$
${\sf reexpand},$ ${\sf expand}$) are only applicable to inputs of the form 
${\sf code}(x) \, 11 \, {\sf code}(y) \, 11 \, z$; so, $S_0$ would end 
before a generator in $\{{\sf decode}_2,$ ${\sf contr},$ ${\sf recontr},$ 
${\sf evR}_{q_2}^{CC},$ ${\sf reexpand},$ ${\sf expand}\}$ can be applied.
Moreover, $S_0$ cannot start with a generator in 
$\{{\sf decode}_2,$ ${\sf contr},$ ${\sf recontr},$ ${\sf evR}_{q_2}^{CC},$
${\sf reexpand},$ ${\sf expand}\}$; indeed, for all inputs 
$\, {\sf code}(w) \, 11 \, u \, \in \, {\sf domC}(X)$, 
$u = {\sf code}(a^m) \, 0010 \,$ contains no 11, so these generators are not 
defined on any element of ${\sf domC}(X)$.
So, $S_0$ is over $\{\pi_0, \pi_1, \rho_0, \rho_1\}$ $\cup$ 
$\{\gamma_j : j \in F\}$.

The actions of $\pi_0, \pi_1, \rho_0, \rho_1$ can change an input in at most
$|S_0|$ positions at the left end of the input, so these actions preserve a 
common suffix $u_2$ and $u$ of length $\ge |u| - |S_0|$.
Thus, if $S_0$ consists only of instances of $\pi_0, \pi_1, \rho_0, \rho_1$, 
the Lemma holds for $S_0$ if $|u| \ge \ell \ge n + |S_0|$ and 
$n > 0$. 

Suppose now that $S_0$ contains $\gamma_j$ for some $j \in F$.  Then (if
$m > |S_0|/|a|$), instances of $\pi_0, \pi_1, \rho_0, \rho_1$ will transform 
the input $\, u = {\sf code}(w) \, 11 \, {\sf code}(a^m) \, 0010 \,$  into a 
word $\, {\sf code}(x) \, 11 \, s \ {\sf code}(a^k) \ 0010 \,$ (for some 
$x, s \in A^*$, $k>0$), such that $\gamma_j$ can be applied. This action 
changes an input in $< |S_0|$ positions at the left end of the input. 
Since $\gamma_j$ is assumed to be applicable now, we must also have $x = j$ 
and $s = y_0 z$ for some $y_0 \in {\sf domC}(\phi_j)$, $z \in A^*$. 
Then the output of $\gamma_j$ is 
$\, \gamma_j\big({\sf code}(j) \, 11 \, s \ {\sf code}(a^k) \ 0010\big)$  
$=$ 
${\sf code}(j) \, 11\, {\sf code}(y_0) \, 11\, z\, {\sf code}(a^k) \, 0010$,
thus the common suffix of $u_2$ and $u$ could decrease by length $\le |y_0|$ 
under the action of $\gamma_j$.  So we let $\ell \ge n + |S_0| + |y_0|$ 
and $n > 0$.  
%Since there are $< 4^{|S_0|}$ possible strings of length $< |S_0|$ over
%$\pi_0, \pi_1, \rho_0, \rho_1$, there are at most that many values for $y_0$
%above; let $k_0$ be the maximum length of
Also, at most one $\gamma_j$ (with $j \in F$) occurs in $S_0$, since after 
$\gamma_j$ the output is of the form  
$\, {\sf code}(w_1) \ 11 \ {\sf code}(u_1) \ 11 \, u_2$, which marks the 
end of the action of $S_0$.   This proves the Lemma for $S_0$.

\smallskip

Inductive step ``$S_i \to S_{i+1}$'', for $0 \leq i < N$:  
By induction we assume that for all $u \in {\sf domC}(\phi_w)$ with 
$|u| \ge \ell$, we have
$\, S_i({\sf code}(w) \, 11 \, u)$ $=$
${\sf code}(w_1) \, 11 \, {\sf code}(u_1) \, 11 \, u_2 \,$ for some 
$w_1, u_1, u_2 \in A^*$, where $u_2$ and $u$ have a common suffix of length 
$\ge n$ ($ > 0$). 
Let us write $S_{i+1} = T_{i+1} S_i$; then $T_{i+1}$ is non-empty (by the
definition of $S_{i+1}$). We also let $T_0 = S_0$.

\smallskip

\noindent {\sf Claim 1:} \ If $T_{i+1}$ contains a generator $g \in$
$\{ {\sf contr},$ ${\sf recontr},$ ${\sf evR}_{q_2}^{CC},$ ${\sf reexpand},$ 
${\sf expand},$ ${\sf decode}_2 \}$, then 
$g$ is the first (i.e., rightmost) letter of $T_{i+1}$, and $g$ occurs only
once. 

Indeed, if $g$ were applicable later in $T_{i+1}$, the output of the 
generator preceding $g$ would be of the form
${\sf code}(w_1) \, 11 \, {\sf code}(u_1) \, 11 \, u_2$, so $S_{i+1}$ would
have ended before $g$ was applied. 

\smallskip

\noindent {\sf Claim 2:} \ If $T_{i+1}$ contains a generator $g \in$ 
$\{ {\sf contr},$ ${\sf recontr},$ ${\sf evR}_{q_2}^{CC},$ ${\sf reexpand},$
${\sf expand}\}$ $\cup$ $\{\gamma_j : j \in F\}$, then 
$g$ is the last (i.e., leftmost) letter of $T_{i+1}$, and $g$ occurs only
once. 

Indeed, such a generator outputs a word of the form 
${\sf code}(w_1) \, 11 \, {\sf code}(u_1) \, 11 \, u_2$. So, $S_{i+1}$ ends 
after such a generator. 

\smallskip

As a consequence of Claims 1 and 2, if $T_{i+1}$ contains a generator 
$g \in$ $\{ {\sf contr},$ ${\sf recontr},$ ${\sf evR}_{q_2}^{CC},$ 
${\sf reexpand},$ ${\sf expand}\}$, 
then $T_{i+1}$ consists of just $g$.
A generator of this form does not change $u_2$.

\smallskip

So we can assume for the remaining cases that $T_{i+1}$ is of the form
$t_{i+1}$, or $\, t_{i+1} \cdot {\sf decode}_2$, or 
$\, \gamma_j \cdot t_{i+1}$, or 
$\, \gamma_j \cdot t_{i+1} \cdot {\sf decode}_2$, where $j \in F$ and  
$t_{i+1}$ is over the generators $\pi_0, \pi_1, \rho_0, \rho_1$.

Let ${\sf code}(w_1) \, 11 \, {\sf code}(u_1) \, 11 \, u_2$ be the input of
$T_{i+1}$ (and this is also the output of $S_i$), where $u_2$ and $u$ have
a common suffix of length $\ge n$.

\smallskip

\noindent $\bullet$
Case where  $T_{i+1}$ is over the generators $\pi_0, \pi_1, \rho_0, \rho_1$:
Then $T_{i+1}$ changes the input in at most $|T_{i+1}|$ positions at the 
left end of the input, so $u_2$ will not be affected if 
$\, \ell - n  \ge |T_{i+1}|$ (and $n > 0$).  

\smallskip

\noindent $\bullet$
Case where $T_{i+1} = t_{i+1} \cdot {\sf decode}_2$, with $t_{i+1}$ over
$\pi_0, \pi_1, \rho_0, \rho_1$:
 The output of ${\sf decode}_2$ is of the form 
$\, {\sf code}(w_1) \, 11 \, u_1 \, u_2$, so the common suffix of $u_2$ and 
$u$ is preserved by ${\sf decode}_2$.
The action of $t_{i+1}$, containing only generators from
$\{\pi_0, \pi_1, \rho_0, \rho_1\}$, affects at most $|t_{i+1}|$ positions 
near the left side of the input, so $u_2$ is not changed if 
$\, \ell - n  \ge |T_{i+1}|$ (and $n > 0$).

\smallskip

\noindent $\bullet$ 
Case where $T_{i+1} = \gamma_j \cdot t_{i+1}$, with $t_{i+1}$ over
$\pi_0, \pi_1, \rho_0, \rho_1$:
Applications of $\pi_0, \pi_1, \rho_0, \rho_1$ change fewer than $|t_{i+1}|$ 
letters of the input near the left end, so the common suffix is not affected 
if $\, \ell - n  \ge |T_{i+1}|$.
When $\gamma_j$ is applied, the output produced will be of the form 
$\, {\sf code}(j)\, 11\, {\sf code}(y_{i+1})\, 11\, z\, {\sf code}(a^n)$
$0010$, where $y_{i+1} \in {\sf domC}(\phi_j)$.  Then $u_2$ will not
be affected if we pick $\ell \ge n + |T_{i+1}| + |y_{i+1}|$ and $n > 0$.

\smallskip

\noindent $\bullet$
Case where $T_{i+1} = \gamma_j \cdot t_{i+1} \cdot {\sf decode}_2$, with 
$t_{i+1}$ over $\pi_0, \pi_1, \rho_0, \rho_1$:
 This case can be handled as a combination of the previous two cases.  
  
\smallskip

\noindent In all the above cases the constraints are fulfilled
for all $i = 0, \ldots, N$, and for all $u = {\sf code}(a^m) \, 0010$, 
if $\, m \ge N + |X| + \sum_{i=0}^N |y_i|$ \, (using the fact that
$\, \sum_{i=0}^N |T_i| = |X|$).  Note that the words $y_i$ do 
not depend on the choice of the input $u = {\sf code}(a^m) \, 0010$,  
whenever $m$ is long enough; indeed, to determine all $y_i$ we can apply
each $S_i$ to the infinite word
$\, {\sf code}(a)^{\omega} \in \{0,1\}^{\omega}$. 
 \ \ \ $\Box$

%\bigskip

\newpage  %%%%%%%%%%%%%%%%%%%%%%%%%%%%%%%%%%%%%%%%%%%%%%%%%%%%

\noindent {\bf Notation.} For a given polynomial $q$ (of the form
$q(n) = a \, n^k + a$ with integers $a, k \ge 1$), let

\medskip

 \ \ \ ${\cal S}_2^{(q)} \ = \ \{f \in {\cal RM}_2^{\sf P} \ : \ f$ 
    is computed by an ${\cal RM}_2^{\sf P}$-program with built-in 
    polynomial $\le q\}$.

\medskip

\noindent We call $w$ an {\em ${\cal S}_2^{(q)}$-program} iff $w$ is an 
${\cal RM}_2^{\sf P}$-program with built-in balance and time-complexity 
polynomial $\le q$.

Let $\, {\cal RM}_2^{(q)} \, = \, \langle {\cal S}_2^{(q)} \rangle$,
 \ i.e., the submonoid of ${\cal RM}_2^{\sf P}$ generated by the set
${\cal S}_2^{(q)}$.
Obviously, we have:

\begin{pro} \label{RMisUnionofRMq} 
 \ For any set of polynomials $\{q_i : i \in {\mathbb N}\}$ of the form
$q_i(n) = a_i \, n^{k_i} + a_i$, such that $\sup\{a_i: i \in {\mathbb N}\}$
$ = +\infty = $  $\sup\{k_i : i \in {\mathbb N}\}$,
we have: 
 \ \ $\bigcup_{i \in {\mathbb N}} \, {\cal RM}_2^{(q_i)} \ = \,$ 
${\cal RM}_2^{\sf P}$.  
 \ \ \ \ \ $\Box$ 
\end{pro} 
The non-finite generation result for ${\cal RM}_2^{\sf P}$ also
holds for ${\cal RM}_2^{(q)}$, and the proof is similar. We need a few
preliminary facts.

\begin{lem} \label{phiVSgamma}
 \ For every polynomial $q$ of the form $q(n) = a \, n^k + a$ with 
$a, k \ge 2$, and every ${\cal S}_2^{(q)}$-program $w$ we have: 
 \ $\gamma_w \in {\cal S}_2^{(q)}$.  
\end{lem} 
{\bf Proof.} Recall that $\gamma_w({\sf code}(w) \, 11 \, uv)$ $=$
${\sf code}(w) \, 11 \, {\sf code}(u) \, 11 \, v$, where 
$u \in {\sf domC}(\phi_w)$.  The input balance of $\gamma_w$ is $\le q$. 
Indeed, the input is shorter than the output; and the output length is 
$\, 2 \, |w| + 2 + 2 \, |u| + 2 + |v|$, which is less than 
$\, q(|{\sf code}(w) \, 11 \, uv|)$ $=$ $q(2 \, |w| + 2 + |u| + 2 + |v|)$ 
when $q(n) \ge 2 \, n^2 + 2$. 

To compute ${\sf code}(w) \, 11 \, {\sf code}(u) \, 11 \, v$ from input 
${\sf code}(w) \, 11 \, uv$, an ${\cal RM}_2^{\sf P}$-machine can proceed 
as follows: 
First, the machine reads and outputs ${\sf code}(w) \, 11$.
Then it runs the program $w$ on input $uv$, i.e., it simulates the
corresponding ${\cal RM}_2^{\sf P}$-machine $M_w$ (which has built-in 
polynomial $q$), with an extra tape and a few modifications. 
While searching for a prefix of $uv$ in ${\sf domC}(\phi_w)$, 
the longest prefix examined so far is kept on the extra tape;
the output $\phi_w(u)$ of $M_w$ will not be written on the output tape.
Once $u$ (the prefix of $uv$ in ${\sf domC}(\phi_w)$) has been found (and
written on the extra tape), ${\sf code}(u) \, 11 \, v$ is appended on the 
output tape. 

All this takes time \ $\le$
$|{\sf code}(w) \, 11| + q(|u|) + |{\sf code}(u) \, 11 \, v|$ $=$ 
$2 \, |w| + 2 + q(|u|) + 2 \, |u| + 2 + |v|$; this is
 \ $< q(2 \, |w| + 2 + |u| + 2 + |v|) = q(|{\sf code}(w) \, 11 \, uv|)$ 
 \ when \ $q(n) \ge 2 \, n^2 + 2$.  
 \ \ \ $\Box$

\begin{lem} \label{GenSetRMq}
 \ Let $q$ be a polynomial that is larger than a certain polynomial of 
degree 5.  Then ${\cal RM}_2^{(q)}$ is generated by  

\hspace{-0.1in} 
$\{\rho_0, \, \rho_1, \, \pi_0, \, \pi_1, \, {\sf decode}_2, $
$ \, {\sf contr}, \,  {\sf recontr}, \, $
${\sf evR}_{q_2}^{CC}, \, {\sf reexpand}, \, {\sf expand}\}$  
$ \ \cup \ $   
$\{\gamma_z: z$ {\rm is an ${\cal S}_2^{(q)}$-program}\}.
\end{lem}
{\bf Proof.} When $w$ is an ${\cal S}_2^{(q)}$-program then as a consequence 
of Lemma \ref{RMsimulation}, 

\medskip

$\phi_w$
$ \ = \ $  $\rho_{{\sf code}(w') \, 11} \circ {\sf decode}_2$ $\circ$
${\sf contr} \circ {\sf recontr}^{2m} \circ {\sf evR}_{q_2}^{CC}$
$\circ$ ${\sf reexpand}^m \circ {\sf expand} \circ \gamma^o_w$,

\smallskip

$ \ = \ $  $\rho_{{\sf code}(w') \, 11} \circ {\sf decode}_2$ $\circ$
${\sf contr} \circ {\sf recontr}^{2m} \circ {\sf evR}_{q_2}^{CC}$
$\circ$ ${\sf reexpand}^m \circ {\sf expand} \circ \gamma_w$ $\circ$ 
$\pi_{{\sf code}(w) \, 11}$,

\medskip

\noindent where $q_2$ is a certain polynomial of degree 2. So the above 
generating set does indeed generate ${\cal RM}_2^{(q)}$. We still need to
show that these generators belong to ${\cal RM}_2^{(q)}$.  

The functions $\rho_0$, $\rho_1$, ${\sf decode}_2$, ${\sf contr}$, 
${\sf recontr}$, ${\sf reexpand}$, ${\sf expand}$, $\pi_0$, $\pi_1$ have 
balance and complexity $\le 4 \, (n+1)^2$.
And $\gamma_w \in {\cal RM}_2^{(q)}$ if $w$ is an ${\cal S}_2^{(q)}$-program 
(by Lemma \ref{phiVSgamma}). 
Let us verify that ${\sf evR}_{q_2}^{CC}$ has balance $\le q_2$ and 
complexity $O(n^5)$. By definition, 

\smallskip
   
${\sf evR}_{q_2}^{CC}({\sf code}(w) \, 11 \, {\sf code}(u) \, 11 \, v)$ 
$=$ ${\sf code}(w) \, 11 \, {\sf code}(\phi_w(u)) \, 11 \, v$. 

\smallskip

\noindent Then ${\sf evR}_{q_2}^{CC}$ has balance $\le q_2$, since on an
output of length $n = 2 \, |w| + 2 + 2 \, |\phi_w(u)| + 2$, the input 
length is $\le$ $2 \, |w| + 2 + 2 \, q_2(|\phi_w(u)|) + 2$ $\le$ 
$q_2\big(2 \, |w| + 2 + 2 \, |\phi_w(u)| + 2 \big) \, = \, q_2(n)$. 

When $\phi_w$ can be computed by an ${\cal RM}_2^{(q)}$-machine with 
built-in polynomial $p_w$ ($\le q_2$), then 
${\sf evR}_{q_2}^{CC}({\sf code}(w) \, 11 \, {\sf code}(u) \, 11)$ can be 
computed in time $\le c \ |w| \ p_w(|u|)^2 \le  c \ |w| \ q_2(|u|)^2$, 
for some constant $c > 0$ (see the proof of Prop.\ 4.4 in \cite{s1f}).  
Since $q_2$ has degree 2, ${\sf evR}_{q_2}^{CC}$ has complexity $O(n^5)$. 
Thus, there exists $q$ of degree 5 such that the above generators belong 
to ${\cal RM}_2^{(q)}$.
 \ \ \ $\Box$

\begin{thm} \label{RMqNotFinGen}
 \ \ For any polynomial $q$ such that $q(n) = a \, n^k + a$, with 
$k \ge 5$ and $a > a_0$ (for some constant $a_0 > 1$), we have: 
 \, ${\cal RM}_2^{(q)}$ is {\em not} finitely generated.
\end{thm}
{\bf Proof.}  The proof is very similar to the proof of Theorem 
\ref{RMnotFinGen}.  We saw in Lemma \ref{GenSetRMq} that 
${\cal RM}_2^{(q)}$ is generated by the infinite set

\smallskip

\hspace{-0.1in} 
$\{\rho_0, \, \rho_1, \, \pi_0, \, \pi_1, \, {\sf decode}_2, $
$ \, {\sf contr}, \,  {\sf recontr}, \, $
${\sf evR}_{q_2}^{CC}, \, {\sf reexpand}, \, {\sf expand}\}$
$ \ \cup \ $
$\{\gamma_z: z$ is an ${\cal S}_2^{(q)}$-program\}.

\smallskip

\noindent Let us assume, by contradiction, that ${\cal RM}_2^{(q)}$ is
finitely generated. Then a finite generating set can be extracted from this
infinite generating set; so ${\cal RM}_2^{(q)}$ is generated by

\smallskip

$\Gamma_{\sf fin} \ = \ $
$\{\rho_0, \, \rho_1, \, \pi_0, \, \pi_1, \, {\sf decode}_2, \, {\sf contr},$
$ \,  {\sf recontr}, \, {\sf evR}_{q_2}^{CC}, \, {\sf reexpand}, $
$ \, {\sf expand} \}$ $ \ \cup \ $ $\{\gamma_i: i \in F\}$,

\smallskip

\noindent where $F$ is some finite set of ${\cal S}_2^{(q)}$-programs.
For every ${\cal S}_2^{(q)}$-program $w$ let $X$ be a word in 
$\Gamma_{\sf fin}^*$ that expresses $\gamma_w$ as a finite sequence of 
generators. 

  From here on, the proof is identical to the proof of Theorem
\ref{RMnotFinGen}.  We use the fact that $\, {\sf domC}(\phi_w)$ $=$
$\{{\sf code}(a^n) \ 0010 \, : \, n > 0\} \,$ is a finite-state language, so
for such a program $w$, $\gamma_w$ has linear complexity (being computable 
by a Mealy machine) and belongs to ${\cal S}_2^{(q)}$.
  \ \ \ $\Box$

%%%%%%%%%%%%%%%%%%%%%%%%%%%%%%%%%%%%%%%%%%%%%%%%%%%%%%%%%%%%%%%%%%%
%%% Section
%%%%%%%%%%%%%%%%%%%%%%%%%%%%%%%%%%%%%%%%%%%%%%%%%%%%%%%%%%%%%%%%%%%

\section{Some complexity consequences of non-finite generation }

%%%%%%%%%%%%%%%%%%%%%%%%%%%%%%%%%%%%%%%%%%%%%%%%%
\subsection{Hierarchy and separation}
%%%%%%%%%%%%%%%%%%%%%%%%%%%%%%%%%%%%%%%%%%%%%%%%%

\begin{pro} \label{fingenRM_q}
 \ Let $q$ be a polynomial of the form $q(n) = a \, n^k + a$ such that
$a, k \ge 1$.
The set ${\cal S}_2^{(q)}$, and hence the monoid ${\cal RM}_2^{(q)}$,
are contained in a finitely generated submonoid of ${\cal RM}_2^{\sf P}$. 
\end{pro}
{\bf Proof.} Let $w$ be a ${\cal RM}_2^{\sf P}$-program such that $\phi_w$ 
has I/O-balance and time-complexity $\le q$. Then ${\sf evR}_q^C$ can 
simulate $\phi_w$ directly, without any need of padding and unpadding. So 
we have for all $u \in {\sf domC}(\phi_w)$, $v \in A^*$:

\smallskip

\hspace{.6in} $\phi_w (u v)$ $ \ = \ $
$\rho_{_{{\sf code}(w) \, 11}} \circ \, {\sf evR}_q^C \, \circ $
$\pi_{_{ {\sf code}(w) \, 11}}(u v)$.

\smallskip

\noindent So ${\cal S}_2^{(q)}$ is contained in the submonoid generated by
$\{\pi_0, \pi_1, \rho_0, \rho_1, {\sf evR}_q^C\}$.  (Compare with Lemma 
\ref{infGenSets} and the proof of Prop.\ 4.5 in \cite{s1f}.)
 \ \ \ $\Box$

\bigskip

\noindent The proof of Prop.\ \ref{fingenRM_q} yields the following chain 
of submonoids in which non-finitely generated and finitely generated
submonoids alternate.

\begin{cor} \label{ChainInRM}
 \ Let \ $\ldots < q_i < q_{i+1} < \ldots$ \ be any sequence of polynomials
such that for all $i \ge 0$, $q_{i+1}$ is is large enough so that
${\sf evR}_{q_i}^C$ has an ${\cal RM}_2^{\sf P}$-program with built-in
polynomial $q_{i+1}$. Then ${\cal RM}_2^{\sf P}$ contains a strict 
inclusion chain, which is infinite in the upward direction,

\medskip

\hspace{0.3in}
$\ldots \ \ \subsetneqq \ {\cal RM}_2^{(q_i)} \ \subsetneqq \ $
$\langle \pi_0, \pi_1, \rho_0, \rho_1,$
        ${\sf evR}_{q_i}^C \rangle_{{\cal RM}_2^{\sf P}} \ $
$\subsetneqq \ {\cal RM}_2^{(q_{i+1})} \ \subsetneqq \ \ \ldots \ \ $
$\ldots$ \ .
\end{cor}
{\bf Proof.} The strictness of the inclusions in the chain follows from
the fact that non-finite generation and finite generation alternate.
  \ \ \ $\Box$

\begin{thm} \label{RMqFinGen} 
 \ Let $q$ be a polynomial of the form $q(n) = a \, n^k + a$ such that 
$a > 1, k \ge 1$.
The submonoid ${\cal RM}_2^{(q)} \subseteq {\cal RM}_2^{\sf P}$ has the 
following properties: \ 

\smallskip

\noindent (1) \ \ ${\cal RM}_2^{(q)}  \neq  {\cal RM}_2^{\sf P}$. 

\smallskip

\noindent (2) \ \ If $q(n) \ge 2 \, (n+1)^2$ (for all $n \in {\mathbb N}$), 
then ${\cal RM}_2^{(q)}$ contains elements of arbitrarily high polynomial 
balance and time-complexity.

\smallskip

\noindent (3) \ \ ${\cal S}_2^{(q)}  \neq  {\cal RM}_2^{(q)}$, if 
$q(n) \ge 2 \, (n+1)^2$.

\medskip

\noindent (4) \ \ ${\sf evR}_q^C \not\in {\cal RM}_2^{(q)}$, if $k \ge 5$
and $a \ge a_0$ (where $a_0$ is as in Theorem \ref{RMqNotFinGen}). 

\smallskip

 \ \ Moreover, ${\sf evR}_q^C$ has balance $\le q$, but its 
time-complexity is {\em not} $\le q$.

\smallskip

\noindent (5) \ \ Let $q_1, q_2$ be polynomials of the above form, such that
$q_1(n) < q_2(n)$ for all $n \in {\mathbb N}$. Suppose also that 
$q_1(n) = a \, n^k + a$ with  $k \ge 5$ and $a \ge a_0$ (as in (4)), 
and that $q_2$ is large enough so that 
${\sf evR}_{q_1}^C \in {\cal RM}_2^{(q_2)}$.
Then $\, {\cal RM}_2^{(q_1)} \subsetneqq {\cal RM}_2^{(q_2)}$.  
\end{thm}
{\bf Proof.} (1) Since ${\cal RM}_2^{(q)}$ is contained in a finitely 
generated submonoid of ${\cal RM}_2^{\sf P}$ (Prop.\ \ref{fingenRM_q}), and 
${\cal RM}_2^{\sf P}$ is not contained (by $\subseteq$) in a finitely 
generated submonoid of ${\cal RM}_2^{\sf P}$ (itself), inequality follows. 

(2) Consider the function $s: 0^n 1 x \mapsto 0^{2 n^2} 1 x$, for all 
$n \ge 0$, $x \in \{0,1\}^*$.  Then $s$ has time-complexity 
$\le 2 \, (n+1)^2$. Indeed, a Turing machine on input $0^n1$ can read this
word $n$ times, each time turning an input $0$ into some new letter $a$,
and each time writing $0^n$ on the output tape; this produces $0^{n^2}$ in 
the output; then one more copy of $0^{n^2}$ is made, followed by $1$. 
This takes time $\le 2 \, (n+1)^2$.  

Then, $s^m$ (i.e., the composition of $m$ instances of $s$) has complexity 
$\ge 2^m \, n^{2^m}$ (since the output length is that 
high, the time must be at least that much too). Thus the functions $s^m$
$\in {\cal RM}_2^{(q)}$ (as $m$ grows) have unbounded complexity, both in 
their degree and in their coefficient.

(3) By (2),  ${\cal RM}_2^{(q)}$ contains functions with arbitrarily high 
polynomial balance and time-complexity, whereas ${\cal S}_2^{(q)}$ only 
contains functions with balance and complexity $\le q$.

(4) By Prop.\ \ref{fingenRM_q}, ${\cal RM}_2^{(q)}$ is contained in the
submonoid generated by $\{\pi_0, \pi_1, \rho_0, \rho_1, {\sf evR}_q^C\}$; 
and we easily see that $\pi_0, \pi_1, \rho_0, \rho_1 \in {\cal RM}_2^{(q)}$.  
Hence, if ${\sf evR}_q^C$ belonged to ${\cal RM}_2^{(q)}$, the monoid 
${\cal RM}_2^{(q)}$ would be finitely generated, contradicting Theorem 
\ref{RMqNotFinGen}.

The input balance of ${\sf evR}_q^C$ is $\le q$ (see Lemma \ref{GenSetRMq},
where this is proved for $q_2$).
It follows that the time-complexity of ${\sf evR}_q^C$ is not $\le q$,
otherwise we would have ${\sf evR}_q^C \in {\cal RM}_2^{(q)}$.  

(5) This follows from (4) since 
$\, {\sf evR}_q^C \not\in {\cal RM}_2^{(q_1)}$, but 
$\, {\sf evR}_q^C \in {\cal RM}_2^{(q_2)}$. 
 \ \ \ $\Box$

%\bigskip

%\noindent {\bf Remark.} Concerning (5) in Theorem \ref{RMqFinGen}: 
%In the proof of Prop.\ 4.4 in \cite{s1f} we saw that if 
%$q_2(n) > c \, n \, (q_1(n))^2$ \, (for a certain constant $c >0$), then 
%${\sf evR}_{q_1}^C \in {\cal RM}_2^{(q_2)}$.
%This gives us a bound on how much is ``large enough'' for $q_2$.  

\begin{cor} \label{RMhierarchy} {\bf (Strict complexity hierarchy of 
submonoids in ${\cal RM}_2^{\sf P}$).} 
 \  There exists an infinite sequence of polynomials 
$(q_i :  i \in {\mathbb N})$, each of the form
$q_i(n) = a_i \, n^{k_i} + a_i$ with $k_i, a_i > 1$, and with
$q_i(n) < q_{i+1}(n)$ for all $i, n \in {\mathbb N}$, such that the 
following holds: 

\medskip

\hspace{0.3in}
${\cal RM}_2^{(q_i)} \subsetneqq {\cal RM}_2^{(q_{i+1})}$ for all $i$, 
and  \ \ \ $\bigcup_{i \in {\mathbb N}} {\cal RM}_2^{(q_i)}$ $=$ 
 ${\cal RM}_2^{\sf P}$.

\medskip

\noindent 
Moreover, ${\cal RM}_2^{\sf P}$ (which is not finitely generated)
is the union of a $\subset$-chain of 4-generated submonoids.
\end{cor}
{\bf Proof.} The first statements follow from Theorem \ref{RMqFinGen} (1) 
and Prop.\ \ref{RMisUnionofRMq}.  The last statement follows from Cor.\
\ref{ChainInRM}.
 \ \ \ $\Box$

%Language complexity hierarchy (if $T_i$ are constructible):
%
%If $(\forall n \in {\mathbb N})[ \, T_1(n) \le T_2(n) \, ]$, and
%$T_2$ is not $O(T_1 \cdot \log T_1)$, then
%${\sf DTime}(T_1) \subsetneqq {\sf DTime}(T_2)$.

\bigskip

Since each ${\cal RM}_2^{(q_i)}$ contains functions of arbitrarily high
polynomial complexity (by Theorem \ref{RMqFinGen} (2)), the monoids 
${\cal RM}_2^{(q_i)}$ form a strict complexity hierarchy of a new sort,
different from the usual complexity hierarchies.  
The fact that ${\cal S}_2^{(q)} \neq {\cal RM}_2^{\sf P}$ could have been 
shown by a diagonal argument. It is not clear whether classical separation 
techniques from complexity theory would show the results (1), (3), (4), (5)
of Theorem \ref{RMqFinGen}.  

\bigskip

\noindent {\bf Remark:}
The monoid {\sf fP}, being finitely generated, does not contain an infinite 
strict complexity hierarchy of monoids (but it can contain hierarchies of 
sets). Indeed, we have in general:

\smallskip

\noindent {\bf Fact.}
{\em A finitely generated monoid $M$ does not contain any infinite strict 
$\omega$-chain of submonoids whose union is $M$.}

\smallskip

\noindent Indeed, if we had a chain $(M_i : i \in \omega)$ with
$M_0 \subsetneqq \ \ldots \ \subsetneqq M_i \subsetneqq M_{i+1} \subsetneqq$
$\ldots \ \ldots \ \subsetneqq \bigcup_{i \in \omega} M_i = M$, then there
would exist $j$ such that $M_j$ contains a finite set of generators of $M$
(since $\bigcup_{i \in \omega} M_i = M$).  Then $M_j = M$, contradicting the 
strict hierarchy. 

\medskip

This Fact does not hold for chains over arbitrary order types; it holds 
for limit ordinals.
The non-finitely generated monoid ${\cal RM}_2^{\sf P}$ contains the 
encoding ${\sf fP}^C$ as a submonoid  (see Section 3 in \cite{s1f}). And 
${\sf fP}^C$ is finitely generated (being an isomorphic copy of {\sf fP}), 
and ${\sf fP}^C$ contains an isomorphic copy of ${\cal RM}_2^{\sf P}$. 
This leads to non-$\omega$ strict chains of submonoids of {\sf fP} and of
${\cal RM}_2^{\sf P}$.

%\bigskip

%\noindent {\bf Remark:} Let $Q_k$ be the set of all polynomials of the
%form $q(n) = a \, n^k + a$ with $a > 1$, for a fixed $k \ge 1$.
%It remains an open the question whether the monoid
% \, $\bigcup_{q \in Q_k} {\cal RM}_2^{(q)}$ \,
%is contained in a finitely generated submonoid of ${\cal RM}_2^{\sf P}$,
%and whether it is finitely generated.

%%%%%%%%%%%%%%%%%%%%%%%%%%%%%%%%%%%%%%%%%%%%%%%%%
\subsection{Irreducible functions} %%%%%%%%%%%%%%%%%%%%%%%%%%%%%%%%%%%%%%%%%%%%%%%%%

Another consequence of non-finite generation is that ${\cal RM}_2^{\sf P}$ 
and ${\cal RM}_2^{(q)}$ have ``irreducible'' elements, i.e., elements that 
cannot be expressed by composition of lower-complexity elements. 
We make this precise in the next definitions.

In this subsection we do not use evaluation maps, so we can use 
``polynomials'' $q(n) = a \, n^k + a$ where we drop the requirement that
$a, k$ are integers, i.e., we now allow real numbers $\ge 1$.

\begin{defn} \label{DEFinfComplexity}
 \ The {\em inf complexity degree} of $f \in {\cal RM}_2^{\sf P}$ is 

\medskip 

$d_f \ = \ $
${\sf inf}\{ k \in {\mathbb R}_{\ge 1} \, : \, f \in {\cal S}_2^{(q)}$
  for some polynomial $q$ of the form $q(n) = b \, n^k + b$,

\hspace{1.4in} for some $b>1$ \}.

\medskip

\noindent We also define the {\em inf complexity coefficient} $c_f$ of $f$
by

\medskip

$c_f \ = \ $
${\sf inf}\{ C_f(\varepsilon) \ : \ \varepsilon \in {\mathbb R}_{> 0} \}$,
 \ \ where

\medskip

$C_f(\varepsilon) \ = \ $
${\sf inf}\{ a \in {\mathbb R}_{\ge 1} \, : \, f \in {\cal S}_2^{(q)}$
  for some polynomial $q$ of the form
  $q(n) = a \, n^{d_f + \varepsilon} + a \}$.

\medskip

\noindent The {\em inf complexity polynomial} of $f$ is the polynomial $q_f$
given by $q_f(n) = c_f \cdot (n^{d_f} + 1)$ (for all $n \in {\mathbb N}$).
\end{defn}
Since $d_f$ and $c_f$ are defined by infimum, $f$ might not be 
in ${\cal S}_2^{(q_f)}$. By the definition of inf we have the following.
%$\varepsilon_1 > 0, \, \varepsilon_2 > 0$:
% \ $f \in {\cal S}_2^{(p_0)}$ where 
% $p_0(n) = (c_f + \varepsilon_2) \cdot (n^{d_f + \varepsilon_1} + 1)$.
% Indeed, by the definition of $C_f(\varepsilon)$,   
% $(C_f(\varepsilon_2) + \frac{1}{2} \varepsilon_1)$
% $(n^{d_f + \varepsilon_2} + 1)$ is a complexity upper-bound for $f$. 
% Moreover, $C_f(\varepsilon_2) < c_f + \frac{1}{2} \varepsilon_1$. 

\begin{pro} 
 \ For any polynomial $q(n) = a \, n^k + a$ with $k > d_f$ and $a > c_f$: 
 \, $f \in {\cal RM}_2^{(q)}$.
 \ \ \ $\Box$
\end{pro}
On the other hand, for every $\varepsilon_1 > 0, \, \varepsilon_2 > 0$:

\smallskip

$f \not\in {\cal S}_2^{(p_1)}$ for any polynomial
$p_1(n) = b \, n^{d_f - \varepsilon_1} + b$ with any $b>1$; 

\smallskip

$f \not\in {\cal S}_2^{(p_2)}$ where
$p_2(n) =  (c_f - \varepsilon_2) \cdot (n^{d_f} + 1)$.

\begin{defn} \label{deltaPrime} 
 \ Let us choose $\delta_1, \delta_2 \in {\mathbb R}_{>0}$.
% $\delta$ with $0 < \delta < \frac{1}{2}$.
A function $f \in {\cal RM}_2^{\sf P}$ is called 
{\em $(\delta_1, \delta_2)$-reducible} iff $f \in {\cal RM}_2^{(q)}$ for some 
polynomial $q(n) = (c_f - \delta_2) \cdot (n^{d_f - \delta_1} + 1)$.  
And $f$ is called {\em $(\delta_1, \delta_2)$-irreducible} iff $f$ is not 
$(\delta_1, \delta_2)$-reducible.
\end{defn}
In other words, $f$ is $(\delta_1, \delta_2)$-reducible iff $f$ is a 
{\em composite} of elements of ${\cal S}_2^{(q)}$ i.e., 
$f \in {\cal RM}_2^{(q)}$, where
$q(n) = (c_f - \delta_2) \cdot (n^{d_f - \delta_1} + 1)$.
So, $f$ can be factored into functions that ``have strictly lower complexity
than $f$'' (regarding both the degree and the coefficient). 
Note that in the definition of $d_f$ and $c_f$ we used ${\cal S}_2^{(q)}$,
not ${\cal RM}_2^{(q)}$ (Def.\ \ref{DEFinfComplexity}).

\begin{pro} \label{ExistsdeltaIrredq1q2} 
 \ For all $\delta_1, \delta_2 \in {\mathbb R}_{>0}$ and all polynomials
$q_1, q_2$ such that ${\cal RM}_2^{(q_1)} \subsetneqq {\cal RM}_2^{(q_2)}$,
there exist $(\delta_1, \delta_2)${\em -irreducible} functions in
 \ ${\cal RM}_2^{(q_2)} - {\cal RM}_2^{(q_1)}$.
\end{pro}
{\bf Proof.} By contradiction, assume that there exist $\delta_1, \delta_2$
such that every $f \in {\cal RM}_2^{(q_2)} - {\cal RM}_2^{(q_1)}$ is 
$(\delta_1, \delta_2)$-reducible, i.e., $f$ can be factored as 
$f = f_m \circ \ \ldots \ \circ f_1$, where  
$f_i \in {\cal RM}_2^{(q_2)}$ ($i = 1, \ldots, m$) with inf degree 
$d_{f_i} < d_f - \delta_1$ and inf coefficient $c_{f_i} < c_f -\delta_2$. 
By the contradiction assumption, among these factors, those that are in 
${\cal RM}_2^{(q_2)} - {\cal RM}_2^{(q_1)}$ can themselves be factored into 
elements of degree and coefficient lower by amount $\delta_1$, respectively 
$\delta_2$. I.e., a factor $f_i \in $ 
${\cal RM}_2^{(q_2)} - {\cal RM}_2^{(q_1)}$ can be factored as 
$f_i = f_{i,m_i} \circ \ \ldots \ f_{i,1}$ with
$d_{f_{i,j}} < d_{f_i} - \delta_1$ and $c_{f_{i,j}} < c_{f_i} - \delta_2$;
hence, $d_{f_{i,j}} < d_f - 2 \, \delta_1$ and
$c_{f_{i,j}} < c_f - 2 \, \delta_2$, for $j = 1, \ldots, m_i$. By repeating 
this process we keep reducing the degree and the coefficient by at least 
$\delta_1$, respectively $\delta_2$, in each step. After a finite number of 
steps we obtain a factorization of $f$ into functions in 
${\cal RM}_2^{(q_1)}$, contradicting the assumption that 
${\cal RM}_2^{(q_1)} \subsetneqq {\cal RM}_2^{(q_2)}$.
 \ \ \ $\Box$

\medskip

\noindent {\bf Remark:} 
A finitely generated monoid, like {\sf fP}, does not contain irreducible 
functions of arbitrarily large complexity. Indeed, all elements are
expressible as a composite of elements of bounded complexity (namely 
the maximum complexity of the finitely many generators).

\bigskip

\medskip

\noindent {\bf Acknowledgement:} The paper benefitted from the referee's
thoughtful reading and advice.

%%%%%%%%%%%%%%%%%%%%%%%%%%%%%%%%%%%%%%%%%%%%%%%%%%%%%%%%%%%%%%%%
%% \small

%%%%%%%%%%%%%%%%%%%%%%%%%%%%%%%%%%%%%%%%%%%%%%%%%%%%

\end{document}